\documentclass[11pt]{article}
\usepackage{benstyle}

\usepackage{overpic}
\usepackage{tikz}
\usepackage{pgfplots}
\usepackage{subfigure}

\newlength{\figurewidth}
\newlength{\figureheight}

\title{
Frames and numerical approximation II: \\
generalized sampling
}
\author{Ben Adcock\\ Department of Mathematics \\ Simon Fraser University \\ Canada \and Daan Huybrechs\\ Department of Computer Science \\ KU Leuven \\ Belgium}
\begin{document}

\maketitle

\begin{abstract}
In a previous paper \cite{FramesPart1} we described the numerical approximation of functions using redundant sets and frames. Redundancy in the function representation offers enormous flexibility compared to using a basis, but ill-conditioning often prevents the numerical computation of best approximations. We showed that, in spite of said ill-conditioning, approximations with regularization may still provide accuracy up to order $\sqrt{\epsilon}$, where $\epsilon$ is a small truncation threshold. When using frames, i.e. complete systems that are generally redundant but which provide infinite representations with coefficients of bounded norm, this accuracy can actually be achieved for all functions in a space. Here, we generalize that setting in two ways. We assume information or samples from $f$ from a wide class of linear operators acting on $f$, rather than inner products associated with the best approximation projection. This enables the analysis of fully discrete approximations based, for instance, on function values only. Next, we allow oversampling, leading to least-squares approximations. We show that this leads to much improved accuracy on the order of $\epsilon$ rather than $\sqrt{\epsilon}$. Overall, we demonstrate that numerical function approximation using redundant representations may lead to highly accurate approximations in spite of having to solve ill-conditioned systems of equations.
\end{abstract}

\section{Introduction}\label{s:introduction}

The approximation of functions in a Hilbert space typically assumes a basis for that space. The non-redundancy of a basis ensures that linear systems associated with approximation problems are non-singular. In addition, suitable structure -- ideally the basis is orthonormal, more generally it may be a Riesz basis -- renders these systems well-conditioned. There is a unique solution, it is stably computable, and there is a close correspondence between properties of the continuous function and of the coefficients in the representation, for example the Parseval identity.

Instead, for a redundant set of functions the corresponding linear systems may be ill-conditioned or even singular, and uniqueness may be lost. Still, good approximations may exist in the span of the set. It may even be much easier to ensure that this is the case than it is for a basis, and in fact this is a popular approach in a wide range of applications. For example, a basis can be `enriched' by adding a few functions that capture a singularity \cite{fix1973singular}. A periodic Fourier basis can be augmented with a few polynomials to capture the non-periodicity of $f$ \cite{roache1978fourier}. Ill-conditioning and redundancy also frequently appear in solution methods for partial differential equations (PDEs).  In Trefftz methods, solutions of a PDE are approximated using other solutions of the same PDE, and this is often successful yet notoriously ill-conditioned \cite{huybrechs2019wbm}. Several methods are based on embedding a domain with complicated geometry $\Omega$, for which a basis is unknown, into a simple bounding box $D$. A basis for $L^2(D)$ yields a (redundant) frame for $L^2(\Omega)$. Examples include embedded/fictitious domain methods, immersed boundary methods and others \cite{boffi2015immersed,kasolis2015fictitious,shirokoff2015volumepenalty}.

As mentioned, redundant representations necessarily lead to ill-conditioning. To which extent are the corresponding function approximations computable? What convergence behaviour and accuracy can be expected? In this paper, we continue a line of investigation that commenced in \cite{FramesPart1} on numerical approximation of functions using redundant function sets in general and frames in particular. The main contribution of \cite{FramesPart1} was a detailed analysis of the accuracy and conditioning of the computation of best approximations with regularization, with a chosen threshold $\epsilon$.

We briefly recall the main results of \cite{FramesPart1} in \S\ref{subsect:bestapproximation}, followed by an overview of the theoretical results of this paper in \S\ref{subsect:results}.

\subsection{Best approximation with regularization}
\label{subsect:bestapproximation}

The main concern of \cite{FramesPart1} was the computation of the best approximation, i.e.\ the orthogonal projection, in the space $\rH_N : = \spn(\Phi_N)$ spanned by a set of $N$ elements $\Phi_N : = \{ \phi_n \}^{N}_{n=1}$.  This approximation is given by $\cP_N f = \sum^{N}_{n=1} x_n \phi_n$, where $\bm{x} = (x_n)^{N}_{n=1}$ is a solution of the linear system
\begin{equation}\label{eq:frame_system}
 \bm{G}_N \bm{x} = \bm{y},
\end{equation}
where $\bm{y} = \{ \ip{f}{\phi_n} \}^{N}_{n=1}$ and $\bm{G}_N = \left \{ \ip{\phi_n}{\phi_m} \right \}^{N}_{m,n=1}$ is the Gram matrix of $\Phi_N$.

If the elements of $\Phi_N$ are nearly or exactly linearly dependent, then $\bm{G}_N$ is ill-conditioned or even singular. Moreover, the coefficients $\bm{x}$ can also grow arbitrarily large, making them impossible to compute in floating point arithmetic for sufficiently large $N$.  The remedy proposed in \cite{FramesPart1} was to regularize \R{eq:frame_system} by using a truncated Singular Value Decomposition (SVD) of $\bm{G}_N$ with a threshold parameter $\epsilon >0$ below which all the singular values are discarded.  This results in a new projection $\cP^{\epsilon}_N f = \sum^{N}_{n=1} (\bm{x}^\epsilon)_n \phi_n$, where $\bm{x}^{\epsilon}$ is the regularized solution of \R{eq:frame_system}.

The main result of \cite{FramesPart1} concerns the best approximation with regularization as follows:
\thm{[\cite{FramesPart1}]
\label{t:Proj_err}
The truncated SVD projection $\cP^{\epsilon}_N$ satisfies
\be{
\label{Proj_err}
\| f - \cP^{\epsilon}_N f \| \leq \inf_{\bm{z} \in \bbC^N} \{ \| f - \sum^{N}_{n=1} z_n \phi_n \| + \sqrt{\epsilon} \| \bm{z} \| \}.
}
Moreover, the (absolute) condition number of the mapping $\bm{y} \mapsto \cP^{\epsilon}_N f$ is at most $1/\sqrt{\epsilon}$.
}
Observe that the right hand side in \eqref{Proj_err} contains two terms. Theorem \ref{t:Proj_err} states that the regularized projection behaves like the best approximation to $f$ in the span of $\Phi_N$ (the first term), as long as the coefficients have sufficiently small norm (second term). Furthermore, convergence can only be expected to an accuracy up to the order of $\sqrt{\epsilon}$. Whether or not this accuracy is achieved, depends on the existence of a representation $\sum^{N}_{n=1} z_n \phi_n$ in the span of $\Phi_N$ with that accuracy and with small norm $\nm{\bm{z}}$ of its coefficients. This question can be studied on a case-by-case basis, as done in \cite{FramesPart1} for a variety of examples.

To answer this question more generally, it is natural to impose additional structure on $\Phi_N$. In  \cite{FramesPart1}, this was done via frames. Recall that an indexed family $\Phi := \{ \phi_n \}_{n=1}^\infty$ is a frame for a Hilbert space $\rH$ if it satisfies the \emph{frame condition}
\begin{equation}
\label{eq:framecondition}
A \nm{f}^2 \leq \sum_{n=1}^\infty | \ip{f}{\phi_n} |^2 \leq B \nm{f}^2, \qquad \forall f \in \rH,
\end{equation}
where $A,B>0$ are positive constants and $\nm{\cdot}$ is the norm on $\rH$. The frame condition ensures the existence of bounded representations, to any accuracy, for all functions in the space. This yields:
\cor{
\label{c:frame_projection}
If $\Phi := \{ \phi_n\}_{n=1}^\infty$ is a frame for $\rH$, and $\Phi_N := \{ \phi_n \}_{n=1}^N$, then
\be{
\label{limit_acc_bad}
\limsup_{N\to\infty} \| f - \cP^{\epsilon}_N f \| \leq \sqrt{\frac{\epsilon}{A}}\, \| f \|, \qquad \forall f \in \rH.
}
}

Unlike in the general setting, the frame condition imposes sufficient structure so that accuracy to order $\sqrt{\epsilon}$ is now guaranteed for all functions in $\rH$. For this reason, as well as the fact that frames occur in numerous computational problems, one can think of frames as an ideal setting in which to apply Theorem~\ref{t:Proj_err}. Of course function approximation with redundancy can be successful without a frame property. For example, in the absence of a frame, one may still use Theorem~\ref{t:Proj_err} to show accuracy in a subspace of $\rH$ consisting of functions with bounded-norm coefficient representations. But this raises the matter of whether functions of interest to a given problem belong to this space. In the absence of a frame structure, this question must then be answered on a case-by-case basis.



\subsection{Main results}
\label{subsect:results}

In this paper, we generalize Theorem \ref{t:Proj_err} using \textit{oversampling} and allowing for \textit{generalized samples} (or \textit{indirect data}). In doing so, we not only allow a much broader class of samples, including, for example, pointwise evaluations, we also overcome the $\sqrt{\epsilon}$ bottleneck.

\pbk
\textbf{Approximation from generalized samples.} In \S \ref{sect:indirect}, instead of inner products with the elements $\phi_n$, as was used in Theorem \ref{t:Proj_err} (see \R{eq:frame_system}), the `data' about the function $f$ is now given by bounded linear functionals $\ell_{m,M} : \rG \rightarrow \bbC$, $m=1,\ldots,M$, which may depend on $M$ and which may be only defined on a dense subspace $\rG$ of $\rH$ (e.g.\ in the case of pointwise evaluations when $\rH = \rL^2(\Omega)$ we consider $\rG=\rC(\overline{\Omega})$).  Very much reminiscent of the frame condition \R{eq:framecondition}, the strongest general statements can be made when this data is sufficiently `rich' so as to stably recover $f$, in particular satisfying
\be{
\label{data_rich_intro}
A' \| f \|^2 \leq \liminf_{M \rightarrow \infty} \sum^{M}_{m=1} | \ell_{m,M}(f) |^2 \leq \limsup_{M \rightarrow \infty} \sum^{M}_{m=1} | \ell_{m,M}(f) |^2 \leq B' \| f \|^2,\quad \forall f \in \rG,
}
for constants $A',B'>0$. We refer to this as a \textit{Marcinkiewicz-Zygmund} condition. A key ingredient is to allow oversampling, i.e.\ let $M >N$, and consider the $M \times N$ linear system
\be{
\label{eq:frame_system_indirect}
\bm{G}_{M,N} \bm{x} \approx \bm{y},\qquad \bm{y} = \{\ell_{m,M}(f) \}^{M}_{m=1},
}
where $\bm{G}_{M,N} = \{ \ell_{m,M}(\phi_n) \}^{M,N}_{m,n=1}$.   As we shall subsequently explain, this system generally remains ill-conditioned for large $N$, even when $M \gg N$. Hence we construct an approximation by singular value thresholding. This leads to a regularized approximation  $\cP^{\epsilon}_{M,N} f$ whose coefficients $\bm{x}^{\epsilon}$ are the solution of the SVD-regularized least-squares problem corresponding to \R{eq:frame_system_indirect}. Our main result for this setup is the following:

\thm{
\label{t:thm2intro}
The truncated SVD projection $\cP^{\epsilon}_{M,N} f$ satisfies
\be{
\label{orange}
\| f - \cP^{\epsilon}_{M,N}f \| \leq \inf_{\bm{z} \in \bbC^N} \left\{   \nm{f - \sum^{N}_{n=1} z_n \phi_n } + \kappa^{\epsilon}_{M,N} \nm{f - \sum^{N}_{n=1} z_n \phi_n}_{M} + \epsilon \lambda^{\epsilon}_{M,N} \| \bm{z} \|  \right\},
}
for constants $\kappa^{\epsilon}_{M,N},\lambda^{\epsilon}_{M,N} > 0$.  The (absolute) condition number of the mapping $\bm{y} \mapsto \cP^{\epsilon}_{M,N} f$ is precisely $\kappa^{\epsilon}_{M,N}$.  Moreover, these constants satisfy
\bes{
\kappa^{\epsilon}_{M,N} \leq \frac{\sqrt{B_N}}{\epsilon},\qquad \lambda^{\epsilon}_{M,N} \leq \frac{\sqrt{B_N}}{\epsilon},
}
for all $M$ and $N$, $M \geq N$, where $B_N$ is the Bessel constant of $\Phi_N$ over $\rH_N = \spn(\Phi_N)$. If the sampling functionals satisfy \R{data_rich_intro} then, for fixed $N$,
\bes{
\limsup_{M \rightarrow \infty} \kappa^{\epsilon}_{M,N} \leq \frac{1}{\sqrt{A'}},\qquad \limsup_{M \rightarrow \infty} \lambda^{\epsilon}_{M,N} \leq \frac{1}{\sqrt{A'}}.
}
}

Here, $\nm{g}^2_{M} = \sum^{M}_{m=1} |\ell_{m,M}(g)|^2$ is the discrete semi-norm defined by the data. We recall that $\Phi_N$ is a Bessel sequence, since it is finite, and therefore it has a finite Bessel constant $B_N > 0$, defined as the smallest constant for which $\sum^{N}_{n=1} | \ip{f}{\phi_n} |^2 \leq B_N \nm{f}^2$, $\forall f \in \rH$.


This result is very general, but its main conclusion is the following. Provided $M$ is sufficiently large and the samples satisfy \R{data_rich_intro} then the approximation error depends on $f - \sum^{N}_{n=1} z_n \phi_n$ (measured in some norm) and $\epsilon \nm{\bm{z}}$. The constants $\kappa^{\epsilon}_{M,N}$ and $\lambda^{\epsilon}_{M,N}$ can be large when $M$ is insufficiently large (behaving like $\sqrt{B_N} / \epsilon$ in the worse case), but they can be made arbitrarily close to $1/ \sqrt{A'}$ by taking $M$ large, i.e. by increasing oversampling. Furthermore, on comparison with Theorem \ref{t:Proj_err}, we notice that $\sqrt{\epsilon}$ in the error bound has been replaced by $\epsilon$. Hence, under suitable assumptions on $f$, $\Phi_N$ and $M$, we expect order $\epsilon$ accuracy in the limit, as opposed to order $\sqrt{\epsilon}$ accuracy.

This result raises the question of how large $M$ needs to be in comparison to $N$. We quantify this via the so-called \textit{stable sampling rate}. Unsurprisingly, the question `how large is sufficiently large' depends completely on system $\Phi_N$ and samples $\{\ell_{m,M} \}$, and thus must be analyzed on a case-by-case basis. In \S \ref{s:SSR}, we illustrate an example for which the stable sampling rate is provably linear, i.e.\ there exists a $C\geq 1$ such that setting $M \geq C N$ implies that the constants $\kappa^{\epsilon}_{M,N}$ and $\lambda^{\epsilon}_{M,N}$ are bounded independently of $\epsilon$. Alternatively, this rate can also be computed numerically, as we explain in \S \ref{ss:SSR}.

\pbk 
\textbf{Frame approximation from frame samples.} Theorem \ref{t:thm2intro} applies for arbitrary $\Phi_N$ and any samples satisfying \R{data_rich_intro}. In order to make statements about the limiting behaviour as $M,N \rightarrow \infty$ we need two ingredients. First, a sequence $\bm{z}$ for which $\sum^{N}_{n=1} z_n \phi_n$ converges to $f$ in the Hilbert space norm and for which the coefficient norm $\nm{\bm{z}}$ does not blow up. Second, additional regularity of the samples and/or $f$ so that the $M$-norm can be controlled by a suitable norm in which one also has $\sum^{N}_{n=1} z_n \phi_n \rightarrow f$.

Whether or not such conditions hold could be answered on a case-by-case basis for particular types of sampling and approximation systems. But instead, we now address them in a general scenario where both the sampling functionals and the approximation system $\Phi_N$ are endowed with a frame structure. Specifically, let $\Phi = \{ \phi_n \}^{\infty}_{n=1}$ and $\Psi = \{ \psi_m \}^{\infty}_{m=1}$ be frames for $\rH$ and set $\Phi_N = \{ \phi_n \}^{N}_{N=1}$ and $\ell_{m,M} = \ip{\cdot}{\psi_m}$ for $m = 1,\ldots,M$. Note that \R{data_rich_intro} now automatically holds with $\rG = \rH$, with $A'$ and $B'$ being the frame bounds for $\Psi$.

Our main result in this case is the following:

\thm{
\label{t:thm3intro}
In the above setting, the truncated SVD projection $\cP^{\epsilon}_{M,N} f$ satisfies
\bes{
\| f - \cP^{\epsilon}_{M,N}f \| \leq  \inf_{\bm{z} \in \bbC^N} \left\{ (1+ \sqrt{B'}\kappa^{\epsilon}_{M,N}) \nm{f - \sum^{N}_{n=1} z_n \phi_n } + \epsilon \lambda^{\epsilon}_{M,N} \| \bm{z} \| \right\},
}
for certain constants $\kappa^{\epsilon}_{M,N},\lambda^{\epsilon}_{M,N} > 0$.  The (absolute) condition number of the mapping $\bm{y} \mapsto \cP^{\epsilon}_{M,N} f$ is precisely $\kappa^{\epsilon}_{M,N}$.  These constants satisfy
\bes{
\kappa^{\epsilon}_{M,N} , \lambda^{\epsilon}_{M,N} \leq \left \{ \begin{array}{cc} \sqrt{B}/\epsilon & \Psi \neq \Phi \\ 1/\sqrt{\epsilon} & \Psi = \Phi,
\end{array} \right . ,
}
for all $M$ and $N$, $M \geq N$, and
\bes{
\limsup_{M \rightarrow \infty} \kappa^{\epsilon}_{M,N} \leq \frac{1}{\sqrt{A'}},\qquad \limsup_{M \rightarrow \infty} \lambda^{\epsilon}_{M,N} \leq \frac{1}{\sqrt{A'}}.
}
Moreover, for each $1 < \theta < \infty$ and $N \in \bbN$ there exists a function $\Theta^{\epsilon}(N,\theta)$ such that
\be{
\label{limit_acc_good}
\limsup_{M \geq \Theta^{\epsilon}(N,\theta)} \| f - \cP^{\epsilon}_{M,N}f \| \leq \epsilon \frac{\theta}{\sqrt{A A'}} \nm{f}.
}
}

The quantity $\Theta^{\epsilon}(N,\theta)$ is what we term the stable sampling rate. As noted above, it can be computed numerically.  Observe that in the special case $\Psi = \Phi$, the setting of Theorem \ref{t:Proj_err} (i.e.\ sampling and approximating with the same functions) is restored. However, oversampling according to \R{limit_acc_good} overcomes the $\sqrt{\epsilon}$ bottleneck in \R{limit_acc_bad}, thus improving the limiting accuracy to order $\epsilon$.

\subsection{Relation to other work}


This paper is a continuation of \cite{FramesPart1}, in which the systematic study of numerical frame approximation was commenced.  This study had its origins in earlier work on so-called \textit{Fourier extensions} \cite{FEStability,huybrechs2010fourier}, which are particular frames arising as restrictions of the Fourier basis on a box to a subdomain.

Our use of oversampling here is inspired by earlier work on \textit{generalized sampling} in Hilbert spaces by the first author and Hansen \cite{BAACHShannon,BAACHAccRecov,BAACHOptimality}.  That work considered both sampling and approximation using orthonormal bases and frames, introducing the stable sampling rate as well, but did not address the ill-conditioning issue for approximation in the latter.  Note that the matrices $\bm{G}_{M,N}$ (in the case of Theorem \ref{t:thm3intro} with $\Psi = \Phi$) and $\bm{G}_N$ are so-called \textit{uneven} and \textit{finite} sections respectively of the infinite Gram matrix of the full frame $\Phi$.  Using uneven as opposed to finite sections is a well-known trick in computational spectral theory \cite{hansen2011,Lindner2008,lindner2006}.

For a more in-depth discussion of relations between this work and standard frame theory, we refer to \cite{FramesPart1}.

Our focus in this paper is accuracy and conditioning of the regularized frame approximations.  We do not consider efficiency, i.e.\ computational time, which is very much dependent on the particular frame under consideration.  There are efficient numerical methods for solving \eqref{eq:frame_system_indirect} for certain frames and sampling functionals \cite{coppe2020splines,coppe2020az,LyonFast,matthysen2015fastfe,matthysen2017fastfe2d}. In the absence of a more efficient method, the cost of computing the SVD of an $M \times N$ matrix with $M > N$ scales as ${\mathcal O}(M N^2)$. However, based in particular on \cite{coppe2020az}, some of the examples at the end of this paper can be implemented in ${\mathcal O}(N (\log N)^2)$ operations. We refer to \cite{coppe2020az} for more examples.

\section{Preliminaries}

We first introduce some notation and useful concepts from frame theory.

\subsection{Bases and frames}
For the remainder of this paper, $\rH$ is a separable Hilbert space over the field $\bbC$. We write $\ip{\cdot}{\cdot}$ and $\nm{\cdot}$ for the inner product and norm on $\rH$ respectively. 
Recall that $\Phi$ is a \textit{Riesz basis} if $\spn(\Phi)$ is dense in $\rH$, and there exist constants $A,B>0$ such that 
\be{
\label{RieszParseval1}
A \| \bm{x} \|^2 \leq \nm{\sum_{n \in I} x_n \phi_n }^2 \leq B \| \bm{x} \|^2,\quad \forall \bm{x} = \{ x_n \}_{n \in I} \in \ell^2(I).
}
Here and throughout, $\ell^2(I)$ denotes the space of square-summable sequences indexed over $I$, and $\nm{\cdot}$ denotes its norm, i.e.\ $\nm{\bm{x}} = \sqrt{\sum_{n \in I} |x_n|^2 }$. Throughout the paper, we assume that the constants appearing in \R{framecond} are the optimal constants such that the inequality holds. We note that a Riesz basis is an orthonormal basis if \R{RieszParseval1} holds with $A = B = 1$. In general, we may view \R{RieszParseval1} as a relaxed version of Parseval's identity.

An indexed family $\Phi$ is a \textit{frame} if
\be{
\label{framecond}
A \| f \|^2 \leq \sum_{n \in I} | \ip{f}{\phi_n} |^2 \leq B \| f \|^2,\quad \forall f \in \rH,
}
for positive constants $A,B>0$. A frame is \textit{tight} if $A = B$. We refer to \R{framecond} as the \textit{frame condition}.  Note that \R{framecond} implies that $\Phi$ is dense in $\rH$.  It follows from \R{RieszParseval1} that a Riesz basis is also a frame with the same constants $A,B$ \cite[Prop.\ 3.6.4]{christensen2003introduction}.  But, in general, a frame need not be a Riesz basis.

\subsection{Operators on frames}

Associated to any frame $\Phi$ (and therefore any Riesz basis) is the so-called \textit{synthesis} operator
\bes{
\cT : \ell^2(I) \rightarrow \rH,\quad \bm{y} = \{ y_n \}_{n \in I} \mapsto \sum_{n \in I} y_n \phi_n.
}
Its adjoint, the \textit{analysis} operator, is given by
\bes{
\cT^* : \rH \rightarrow \ell^2(I),\quad f \mapsto \{ \ip{f}{\phi_n} \}_{n \in I},
}
and the composition $\cS = \cT \cT^*$, known as the \textit{frame} operator, is 
\bes{
\cS : \rH \rightarrow \rH,\quad f \mapsto \sum_{n \in I} \ip{f}{\phi_n} \phi_n.
}
This operator is self-adjoint, bounded, invertible and positive with 
\be{
\label{Sbounds}
A \nm{f}^2 \leq \ip{\cS f}{f} \leq B \nm{f}^2.
}
See \cite[Lemma 5.1.5]{christensen2003introduction}.  Note that this inequality is equivalent to the frame condition \R{framecond}.  Note also that $\cS = \cI$ is the identity operator for an orthonormal basis.  Similarly, $\cS = A \cI$ for a tight frame.  However, for a general Riesz basis or frame, $\cS \neq \cI$. 


\subsection{Dual frames}
\label{ss:canonical_dual_frame}

A frame $\Psi = \{ \psi_n \}_{n \in I} \subseteq \rH$ is a \textit{dual frame} for a given frame $\Phi$ if 
\be{
\label{dual_frame}
f = \sum_{n \in I} \ip{f}{\psi_n} \phi_n = \sum_{n \in I} \ip{f}{\phi_n} \psi_n,\qquad \forall f \in \rH.
}
If a frame $\Phi$ is also a Riesz basis then it has a unique dual frame $\Psi$, which is also a Riesz basis.  In this case, the pair $(\Phi,\Psi)$ is \textit{biorthogonal}:
\bes{
\ip{\phi_n}{\psi_m} = \delta_{n,m}, \quad n,m \in I.
}
Note that an orthonormal basis is self-dual, i.e. $\Psi = \Phi$. Conversely, a frame and any of its duals are typically not biorthogonal. A frame may also have more than one dual. The so-called \textit{canonical} dual frame of $\Phi$ is the frame $\Psi = \{ \cS^{-1} \phi_n \}_{n \in I}$. This is a frame \cite[Lem. 5.1.5]{christensen2003introduction}, and its frame bounds are $B^{-1}$ and $A^{-1}$ respectively.  In this case, \R{dual_frame} reads
\be{
\label{dual_rep}
f = \sum_{n \in I} \ip{f}{\cS^{-1} \phi_n} \phi_n = \sum_{n \in I} \ip{\cS^{-1} f}{\phi_n} \phi_n.
}
We refer to the coefficients $\bm{a} = \{ \ip{f}{\cS^{-1} \phi_n } \}_{n \in I}$ as the \textit{canonical frame coefficients} of $f$.  These coefficients have the property that, amongst all possible representations of $f$ in $\Phi$, they have the smallest norm \cite[Lem.\ 5.4.2]{christensen2003introduction}. Specifically, if $f = \sum_{n \in I} a_n \phi_n = \sum_{n \in I} c_n \phi_n$ for some $\bm{c} = \{ c_n \}_{n \in I}$, then $\| \bm{c} \| \geq \| \bm{a} \|$. Moreover, from the frame condition of the dual we have $\| \bm{a} \| \leq \frac{1}{\sqrt{A}} \| f \|$.

\section{Approximation from indirect data}
\label{sect:indirect}

In this section, we describe the general setup, which will lead eventually to Theorem \ref{t:thm2intro}. Throughout this section, the \textit{approximation system} is defined by $\Phi_N = \{ \phi_n \}_{n \in I_N } \subset \rH$. It is an an indexed family of $N$ elements in $\rH$, where $I_N$ is an index set of cardinality $N$. For convenience, we now make a mild generalization over \S \ref{s:introduction}, allowing $I_N$ to be an arbitrary index set rather than just $\{1,\ldots,N\}$. As noted, $\Phi_N$ is a Bessel sequence, since it is finite. We write $B_N>0$ for its Bessel constant, i.e.\ the smallest constant for which $\sum^{N}_{n=1} | \ip{f}{\phi_n} |^2 \leq B_N \nm{f}^2$, $\forall f \in \rH$. We also define the operators
\eas{
&\cT_N :  \bbC^{N} \rightarrow \rH,\ \bm{z} = (z_n)_{n \in I_N} \mapsto \sum_{n \in I_N} z_n \phi_n
\\
&\cT^*_N : \rH \rightarrow \bbC^N,\ f \mapsto \left ( \ip{f}{\phi_n} \right )_{n \in I_N}
\\
&\cS_N = \cT_N \cT^*_N : \rH \rightarrow \rH,\ f \mapsto \sum_{n \in I_N} \ip{f}{\phi_n} \phi_n.
}


\subsection{Indirect data}
\label{subsect:indirect_prelim}

Let $\rG$ be a dense subspace of the Hilbert space $\rH$ endowed with a norm $\tnm{\cdot}$.  Suppose that $f$, the function we seek to approximate, and $\Phi_N$ both belong to $\rG$.  For each $M \in \bbN$ let $J_M$ be an index set of cardinality $|J_M| = M$, and 
\bes{
\ell_{m,M} : \rG \rightarrow \bbC,\quad m \in J_M,
}
be a set of linear functionals which are bounded with respect to $\tnm{\cdot}$, i.e.
\be{
\label{functional_bound}
| \ell_{m,M}(f) | \leq c_{m,M} \tnm{f}, \qquad f \in \rG.
}
The data of $f$ is given by
\bes{
\bm{y} = \{ \ell_{m,M}(f) \}_{m \in J_M}.
}
Write $\cM_{M} : \rG \rightarrow \bbC^M$ for the mapping $\cM_{M} f = \{ \ell_{m,M}(f) \}_{m \in J_M}$.   Our goal is to compute an approximation to $f$ in $\Phi_N$ for some $N \leq M$ from this data.

In order to make meaningful general statements about the subsequent approximations we define, we require the data to be sufficiently rich.  In analogy to the frame bounds \R{framecond}, we shall often assume that there exist constants $A',B' > 0$ such that
\be{
\label{data_rich}
A' \| f \|^2 \leq \liminf_{M \rightarrow \infty} \sum_{m \in J_M} | \ell_{m,M}(f) |^2 \leq \limsup_{M \rightarrow \infty} \sum_{m \in J_M} | \ell_{m,M}(f) |^2 \leq B' \| f \|^2,\quad \forall f \in \rG.
}
We term this a \textit{Marcinkiewicz-Zygmund} condition, because it is similar (although not identical) to the classical Marcinkiewicz-Zygmund inequality in approximation theory \cite{LubinskyMarcinkiewicz}, see also \cite{GrochenigMZineq}.
We comment further on this assumption and the constants involved in \S\ref{ss:indirect_constants}. 

Before going any further, let us mention several examples of this framework:

\examp{
Suppose that $M = N$ and the samples of $f$ are inner products with respect to the functions $\phi_m$, i.e.\ $\ell_{m,M}(f) = \ip{f}{\phi_m}$, $m=1,\ldots,M$.  Then this is precisely the setting of Theorem \ref{t:Proj_err}, and is a special case of the present setup with $\rG = \rH$, $\tnm{\cdot} = \nm{\cdot}$.  Note that \R{data_rich} is precisely the frame condition \R{framecond}; in particular, it holds with $A' = A$ and $B' = B$. 
}

\examp{
\label{ex:otherframe}
Let $\Psi = \{ \psi_m \}_{m \in J}$ be another frame (or Riesz/orthonormal basis) of $\rH$ and consider samples of the form $\ell_{m,M}(f) = \ip{f}{\psi_m}$, $m \in J_M$, where $\{ J_M \}_{M \in \bbN}$ are nested index sets with $|J_M | = M$ and $\bigcup^{\infty}_{M=1} J_M = J$.  This problem corresponds to sampling  according to the frame $\Psi$ and reconstructing in $\Phi$, as in the original framework of generalized sampling \cite{BAACHShannon,BAACHAccRecov}.  We also have $\rG = \rH$ and $\tnm{\cdot} = \nm{\cdot}$ in this case, and \R{data_rich} holds with $A'$, $B'$ being the frame bounds for $\Psi$. In fact, it is straightforward to see that the the upper bound holds for any $M$ in this case, i.e.
\be{
\label{framesamplesUB}
 \| f \|_M \leq \sqrt{B'} \| f \|, \qquad \forall f \in \rH.
}
}

\examp{
\label{ex:equispaced}
Consider a frame for the Hilbert space $\rH = \rL^2(-1,1)$ of square-integrable functions on $(-1,1)$.  Suppose that the samples are pointwise evaluations at equispaced points.  In this case, we may take $\rG = \rC([-1,1])$ with its usual norm, and $\ell_{m,M}(f) = \sqrt{2/M} f(-1+2(m-1)/M)$, $m=1,\ldots,M$.  Then \R{data_rich} holds with $A' = B' = 1$, since $\sum^{M}_{m=1} | \ell_{m,M}(f) |^2$ is a Riemann sum approximation to $\int^{1}_{-1} | f(x) |^2 \D x$. More generally, if $-1 \leq t_{1,M} < \ldots < t_{M,M} \leq 1$ are (not necessarily equispaced) sample points, then \R{data_rich} can be achieved with $A' = B' = 1$ if 
\bes{
h_{M} = \max_{0,\ldots,M} \{ t_{m+1,M} - t_{m,M} \} \rightarrow 0,\qquad M \rightarrow \infty,
}
where $t_{0,M} = t_{M,M}-2$ and $t_{M+1,M} = t_{1,M}+2$. Indeed, in this case choosing the linear functionals as
\bes{
\ell_{m,M}(f) = \sqrt{\frac12(t_{m+1,M} - t_{m-1,M})} f(t_{m,M})
}
gives a convergent approximation $\sum^{M}_{m=1} | \ell_{m,M}(f) |^2$ to $\nm{f}^2$ as $M \rightarrow \infty$.
}

\subsection{Best approximation with regularization from indirect data}

Given $f \in \rG$ and data $\cM_{M}f$, we construct the approximation $\cP^{\epsilon}_{M,N}f$ as follows.  Let
\bes{
\bm{G}_{M,N} = \cM_M \cT_N = \{ \ell_{m,M}(\phi_n) \}_{m \in J_M,n \in I_N} \in \bbC^{M \times N}.
}
As we explain below, much like in the setting of Theorem \ref{t:Proj_err}, this matrix is generally ill-conditioned. This arises from the inherent redundancy of $\Phi_N$, independently of the samples -- in particular, it cannot be avoided by taking $M \gg N$.  Hence we need to regularize.  Suppose that $\bm{G}_{M,N}$ has singular value decomposition 
\be{
\label{eq:Gsvd}
\bm{G}_{M,N} = \bm{U} \bm{\Sigma} \bm{V}^*.
}
Let $\epsilon > 0$ be fixed.  Then we set
\be{
\label{Ldata_def}
\bm{x}^{\epsilon} = (\bm{G}^{\epsilon}_{M,N})^{\dag} \bm{y} = \bm{V} (\bm{\Sigma}^{\epsilon})^{\dag} \bm{U}^* \bm{y},
}
where $\dag$ denotes the pseudoinverse and $\bm{\Sigma}^{\epsilon}$ is the diagonal matrix with $n^{\rth}$ entry $\sigma_n$ if 
\be{
\label{sigma_cutoff}
\sigma_n >  \epsilon,
}
and zero otherwise.
The corresponding approximation to $f$ is
\bes{
f \approx \cP^{\epsilon}_{M,N} f = \cT_N \bm{x}^{\epsilon}.
}
Observe that both the solution vector $\bm{x}^{\epsilon}$ and the approximation $\cP^{\epsilon}_{M,N} f$ are uniquely defined by construction, even if the family $\Phi_N$ happens to be linearly dependent. For convenience, we now define the mapping
\bes{
\cL^{\epsilon}_{M,N} : \bbC^M \rightarrow \rH_N,\ \bm{y} \mapsto \cT_N (\bm{G}^{\epsilon}_{M,N})^\dag \bm{y},
}
so that
\be{
\label{Peps_data_def}
\cP^{\epsilon}_{M,N} f = \cT_{N} \bm{x}^{\epsilon} = \cL^{\epsilon}_{M,N} \cM_M f.
}

\rem{
To see why $\bm{G}_{M,N}$ is generally ill-conditioned, consider the setting of Example \ref{ex:otherframe} and suppose that $\Psi$ is a tight frame, i.e.\ $A' = B'$. This assumption is made for convenience: the following argument readily extends to general frames.
Let $\widetilde{\cS}_M : f \mapsto \sum_{m \in J_M} \ip{f}{\psi_m} \psi_m$ be the partial frame operator for $\Psi$. Since $\Psi$ is tight, $\cS_{M} \rightarrow A'\cI$ strongly as $M \rightarrow \infty$, and therefore, for fixed $N$,
\bes{
(\bm{G}_{M,N})^* \bm{G}_{M,N} \rightarrow A' \bm{G}_N,\qquad M \rightarrow \infty,
}
where $\bm{G}_N$ is the Gram matrix of $\Phi_N$. Hence, whenever $\bm{G}_N$ is ill-conditioned (i.e.\ $\Phi_N$ is near-redundant), we expect $\bm{G}_{M,N}$ to inherit the same ill-conditioning for large $M$.
}

\rem{
This argument gives some insight into the advantage of oversampling.  For a tight frame, $(\bm{G}_{M,N})^* \bm{G}_{M,N}$ is an approximate factorization of $\bm{G}_N$.  Thus, solving the linear system \R{eq:frame_system} is akin to solving the normal equations of the least-squares problem $\bm{G}_{M,N} \bm{x} \approx \bm{y}$.  In this sense it is not surprising that oversampling yields $\ord{\epsilon}$ accuracy, whereas solving \R{eq:frame_system} yields only $\ord{\sqrt{\epsilon}}$ accuracy.  Indeed, this is reminiscent of the typical squaring of the condition number incurred when forming the normal equations of a least-squares problem \cite[\S5.3]{golub1996matrix}.
}

\subsection{The solution as an orthogonal projection}

A key element of our subsequent analysis is the reinterpretation of the operator $\cP^{\epsilon}_{M,N}$ as a projection with respect to a semi-definite sesquilinear form. Specifically, we now define the data-dependent sesquilinear form $\ip{\cdot}{\cdot}_M$ on $\rG \times \rG$ as
\bes{
\ip{f}{g}_M = \ip{\cM_M f}{\cM_M g} = \sum_{m \in J_M} \ell_{m,M}(f) \overline{\ell_{m,M}(g)},\quad f,g \in \rG,
}
with corresponding discrete semi-norm $\nm{f}_M = \sqrt{\ip{f}{f}_M} = \nm{\cM_M f}$. Note that in general $\ip{\cdot}{\cdot}_M$ is semi-definite on $\rG_N \times \rG_N$ as well, since
\bes{
\ip{g}{g}_M = \nm{g}_M^2 = \sum_{m \in J_M} |\ell_{m,M}(g) |^2 \geq 0,\quad \forall g \in \rH_N, \quad g \neq 0.
}
In particular, for poorly chosen functionals it may be that $\Vert g \Vert_M = 0$ for some non-trivial function $g \in \rH_N$. However, with assumption \eqref{data_rich} we do have the limiting behaviour
\bes{
 \liminf_{M \to \infty} \ip{g}{g}_M = \liminf_{M \to \infty} \nm{g}_M^2 \geq A' \| g \|^2 > 0,\quad \forall g \in \rH_N, g \neq 0.
}
This means that, ultimately, the sesquilinear form becomes an inner product on all of $\rH_N$.



Recall the singular value decomposition \eqref{eq:Gsvd}, and let $\bm{u}_1,\ldots,\bm{u}_M$, $\bm{v}_1,\ldots,\bm{v}_N$ and $\sigma_1,\ldots,\sigma_N$ be the left and right singular vectors and singular values of $\bm{G}_{M,N}$ respectively, with $\sigma_n \geq 0$, $n=1,\ldots,N$. To the right singular vectors, we associate the functions
\be{\label{eq:xi}
\xi_n = \cT_N \bm{v}_n \in \rH_N,\quad n=1,\ldots,N.
}
By construction, these functions are orthogonal with respect to $\ip{\cdot}{\cdot}_M$. Indeed, by orthogonality of the singular vectors, we have
\be{
\label{OS_orthog}
\ip{\xi_m}{\xi_n}_M = \ip{\cM_M \cT_N \bm{v}_n}{\cM_M \cT_N \bm{v}_n} = \sigma_n \sigma_m \ip{\bm{u}_m}{\bm{u}_n} = \sigma_n \sigma_m \delta_{m,n},\quad m,n \in I_N.
}
Here, too, it may be that $\| \xi_n \|_M = 0$. This is the case if $\sigma_n = 0$.


We shall, for convenience, let $\bm{x} = \bm{x}^{0}$ be the solution of the unregularized problem, given by
\bes{
\bm{x} = (\bm{G}_{M,N})^{\dag} \bm{y}.
}
We also write $\cP_{M,N} = \cP^{0}_{M,N}$ so that $\cP_{M,N} f = \cT_N \bm{x}$.
Using the expression for the pseudoinverse in terms of the SVD, we can write both $\bm{x}$ and $\bm{x}^{\epsilon}$ in terms of the left and right singular vectors:
\be{
\label{x_SVD_def}
\bm{x} = \sum_{\sigma_n > 0} \frac{\ip{\bm{y}}{\bm{u}_n}}{\sigma_n} \bm{v}_n,\qquad \bm{x}^{\epsilon} = \sum_{\sigma_n > \epsilon} \frac{\ip{\bm{y}}{\bm{u}_n}}{\sigma_n} \bm{v}_n.
}
Furthermore, for $\sigma_n > 0$ we also have
\bes{
\ip{\bm{y}}{\bm{u}_n} = \frac{\ip{\cM_M f}{\bm{G}_{M,N} \bm{v}_n}}{\sigma_n} = \frac{\ip{\cM_M f}{ \cM_M \cT_N \bm{v}_n}}{\sigma_n} = \frac{\ip{f}{\xi_n}_M}{\sigma_n},
}
where in the last step we use \R{eq:xi}.  In particular, this gives
\be{
\label{SVD_proj_def}
\cP^{\epsilon}_{M,N} f = \cT_N \bm{x}^{\epsilon} = \sum_{\sigma_n > \epsilon} \frac{\ip{\bm{y}}{\bm{u}_n}}{\sigma_n} \cT_N \bm{v}_n = \sum_{\sigma_n > \epsilon} \frac{\ip{f}{\xi_n}_M}{\sigma^2_n} \xi_n.
}
Similarly, we have
\bes{
\cP_{M,N} f = \cT_N \bm{x} = \sum_{\sigma_n > 0} \frac{\ip{f}{\xi_n}_M}{\sigma^2_n} \xi_n.
}
Finally, we define the regularized spaces
\bes{
\rH^{\epsilon}_{M,N} = \spn \left \{ \xi_n : \sigma_n > \epsilon \right\}.
}
Since $\{ \bm{v}_n \}_{n \in I_N}$ is a basis of $\bbC^N$, we see that the functions  $\{ \xi_n \}_{\sigma_n > 0}$ are an orthogonal basis of $\rH_{M,N}^0$ with respect to $\ip{\cdot}{\cdot}_M$. We can conclude that $\cP_{M,N}$ is the orthogonal projection onto $\rH_{M,N}^0 \subseteq \rH_N \subset G$. In turn, $\cP^{\epsilon}_{M,N}$ is the orthogonal projection onto the subspace $\rH^{\epsilon}_{M,N} \subseteq \rH^0_{M,N}$.

A relevant property in the analysis that follows, is that these orthogonal projections imply a reduction in the $M$-norm:
\lem{
\label{l:projection}
For any $\epsilon \geq 0$, we have
\be{
\label{Mnorm_reduction}
 \| P_{M,N}^\epsilon f \|_M \leq \| f \|_M, \qquad \forall f \in G.
}
}

\subsection{Theoretical results}\label{ss:indirect_analysis}
We now define the constants
\be{
\label{kappa_lambda_indirect}
\kappa^{\epsilon}_{M,N} = \max_{\substack{\bm{y} \in \mathrm{Ran}(\cM_M) \\ \| \bm{y} \| =1}} \| \cL^{\epsilon}_{M,N} \bm{y} \|,\quad \lambda^{\epsilon}_{M,N} = \epsilon^{-1} \max_{\substack{\bm{z} \in \bbC^N \\ \| \bm{z} \| =1}} \| \cT_N  \bm{z} - \cP^{\epsilon}_{M,N}  \cT_N \bm{z} \|.
}
Note that $\kappa^{\epsilon}_{M,N} $ is precisely the operator norm of $\cL_{M,N}^{\epsilon} : \mathrm{Ran}(\cM_M) \rightarrow \rH_N$.  Since $\cL_{M,N}^{\epsilon}$ is linear, it is also its absolute condition number, i.e.\ $\kappa^{\epsilon}_{M,N} $ measures the absolute effect of perturbations in the data $\bm{y}$ on the final approximation.  The constant $\lambda^{\epsilon}_{M,N}$ measures how close $\cP^{\epsilon}_{M,N}$ is to being a projection on the subspace $\rH_N = \cT_N(\bbC^N)$.  


Our first result concerns the approximation error of $\cP^{\epsilon}_{M,N}f$:
\thm{
\label{t:indirect_err}
Let $f \in \rG$.
The truncated SVD approximation $\cP^{\epsilon}_{M,N}f$ satisfies
\be{
\label{indirect_err}
\| f - \cP^{\epsilon}_{M,N}f \| \leq \inf_{\bm{z} \in \bbC^N} \left\{  \| f - \cT_N \bm{z} \| + \kappa^{\epsilon}_{M,N} \| f - \cT_N \bm{z} \|_M + \epsilon \lambda^{\epsilon}_{M,N} \| \bm{z} \| \right\}.
}
}

This result differs from Theorem \ref{t:Proj_err} in several respects. On the one hand, if the constants $\kappa^{\epsilon}_{M,N}$ and $\lambda^{\epsilon}_{M,N}$ are order one, the error depends on $\epsilon \nm{\bm{z}}$, not $\sqrt{\epsilon} \nm{\bm{z}}$, thus overcoming the $\sqrt{\epsilon}$ bottleneck. We will discuss when this occurs in the next subsection. On the other hand, the error bound involves the discrete data norm $\| f - \cT_N \bm{z} \|_M$. In general, this cannot be bounded by $\| f - \cT_N \bm{z}\|$. However, one clearly has
\bes{
\| f - \cT_N \bm{z} \|_M \leq C_M \tnm{f - \cT_N \bm{z}},\qquad C_M = \sqrt{\sum_{m \in J_M} (c_{m,M})^2},
}
where $c_{m,M}$ are the norms of the functionals $\ell_{m,M}$; recall \eqref{functional_bound}.
In particular, for Example \ref{ex:equispaced} it follows that $\| f - \cT_N \bm{z} \|_M \leq \sqrt{2} \| f - \cT_N \bm{z} \|_{\rL^\infty}$.

\prf{[Proof of Theorem \ref{t:indirect_err}]
For any $\bm{z} \in \bbC^N$,
\bes{
\| f - \cP^{\epsilon}_{M,N}f \| \leq \| f - \cT_N \bm{z} \| + \| \cP^{\epsilon}_{M,N}(f-\cT_N \bm{z}) \| + \| \cT_N \bm{z} - \cP^{\epsilon}_{M,N}  \cT_N \bm{z} \|.
}
Consider the second term.  We have
\bes{
\| \cP^{\epsilon}_{M,N}(f-\cT_N \bm{z}) \| = \|  \cL^{\epsilon}_{M,N}  \cM_M (f-\cT_N \bm{z}) \| \leq \kappa^{\epsilon}_{M,N} \| \cM_M (f-\cT_N \bm{z}) \| = \kappa^{\epsilon}_{M,N} \| f - \cT_N \bm{z} \|_M,
} 
which gives the corresponding term in \R{indirect_err}.  Now consider the third term.  It follows immediately from the definition of $\lambda^{\epsilon}_{M,N}$ that $\| \cT_N \bm{z} - \cP^{\epsilon}_{M,N}  \cT_N \bm{z} \| \leq \epsilon \lambda^{\epsilon}_{M,N} \| \bm{z} \|$, as required.
}

We now consider the coefficients $\bm{x}^{\epsilon}$:

\thm{
\label{t:coefficients_indirect}
Let $f \in \rG$.
The coefficients $\bm{x}^{\epsilon}$ of the truncated SVD projection $\cP^{\epsilon}_{M,N}$ satisfy
\be{
\label{coeffics_bounded_indirect}
\| \bm{x}^{\epsilon} \| \leq \inf_{ \bm{z} \in \bbC^N} \left\{ \frac{1}{\epsilon} \,  \| f - \cT_N \bm{z} \|_M + \| \bm{z} \| \right\}.
}
Moreover, if $\Phi = \{ \phi_n \}_{n \in I}$ is a frame, $\Phi_N = \{ \phi_n \}_{n \in I_N}$ and if $\bm{a} =  \{ \ip{f}{\cS^{-1}\phi_n} \}_{n \in I}$ are the canonical frame coefficients of $f$ and $\bm{a}^{\epsilon}_{M,N} \in \ell^2(I)$ is the extension of $\bm{x}^\epsilon$ by zero, then
\be{
\label{coeffics_conv3}
\| \bm{a} - \bm{a}^{\epsilon}_{M,N} \| \leq \sqrt{\sum_{n \in I \backslash I_N} | a_n |^2} + \frac{1}{\epsilon} \| (\cS - \cS_N) \cS^{-1} f \|_M +  \epsilon \frac{\lambda^{\epsilon}_{M,N}}{\sqrt{A}} \| \bm{a} \|.
}
}

For general measurements, \R{coeffics_conv3} does not imply convergence of the coefficients $\bm{a}^{\epsilon}_{M,N}$ to the canonical frame coefficients $\bm{a}$ since the term $\| (\cS - \cS_N) \cS^{-1} f \|_M$ cannot be bounded by $\| (\cS - \cS_N ) \cS^{-1} f \|$ in general.  There is also no guarantee that $\tnm{(\cS - \cS_N ) \cS^{-1} f } \rightarrow 0$ as $N \rightarrow \infty$.  
This does hold, however, when the data arises from sampling with another frame $\{ \psi_n \}_{n \in I}$, as in Example \ref{ex:otherframe}, due to \R{framesamplesUB}. We will discuss this case further in \S \ref{s:framesforall}.

\prf{[Proof of Theorem \ref{t:coefficients_indirect}]
For the first part, we use \R{x_SVD_def} to write
\be{
\label{x_OS_split}
\bm{x}^{\epsilon} = \sum_{\sigma_n > \epsilon} \frac{\ip{f}{\xi_n}_M}{\sigma^2_n} \bm{v}_n =\sum_{\sigma_n > \epsilon} \frac{\ip{f-\cT_N \bm{z}}{\xi_n}_M}{\sigma^2_n} \bm{v}_n + \sum_{\sigma_n > \epsilon} \frac{\ip{\cT_N \bm{z}}{\xi_n}_M}{\sigma^2_n} \bm{v}_n.
}
Consider the first term on the right-hand side.  Since the $\bm{v}_n$ are orthonormal, we have
\bes{
\nm{\sum_{\sigma_n > \epsilon} \frac{\ip{f-\cT_N \bm{z}}{\xi_m}_M}{\sigma^2_n} \bm{v}_n}^2 = \sum_{\sigma_n > \epsilon} \frac{|\ip{f-\cT_N \bm{z}}{\xi_m}_M|^2}{\sigma^4_m} \leq \frac{1}{\epsilon^2} \sum_{\sigma_n > \epsilon} \frac{|\ip{f-\cT_N \bm{z}}{\xi_m}_M|^2}{\sigma^2_m}.
}
It follows from \R{OS_orthog} and \R{SVD_proj_def} that
\bes{
\sum_{\sigma_n > \epsilon} \frac{|\ip{g}{\xi_m}_M|^2}{\sigma^2_m} = \nm{\cP^{\epsilon}_{M,N} g}^2_{M},\quad g \in \rG.
}
Hence 
\bes{
\nm{\sum_{\sigma_n > \epsilon} \frac{\ip{f-\cT_N \bm{z}}{\xi_m}_M}{\sigma^2_n} \bm{v}_n}^2 \leq \frac{1}{\epsilon^2} \nm{\cP^{\epsilon}_{M,N}(f-\cT_N \bm{z})}_{M} \leq \frac{1}{\epsilon^2} \| f - \cT_N \bm{z} \|^2_M,
}
where in the second step we use the fact that $\cP^{\epsilon}_{M,N}$ is the orthogonal projection with respect to $\ip{\cdot}{\cdot}_M$ (recall Lemma~\ref{l:projection}).  This gives the first term of \R{coeffics_bounded_indirect}.  Next, consider the second term of the right-hand side of \R{x_OS_split}.  Since
\be{
\label{useful1}
\ip{\cT_N \bm{z}}{\xi_m}_M = \ip{\cT_N \bm{z}}{\cT_N \bm{v}_n }_M = \sigma^2_n \ip{\bm{z}}{\bm{v}_n},
}
it follows that
\bes{
\nm{\sum_{\sigma_n > \epsilon} \frac{\ip{\cT_N \bm{z}}{\xi_m}_M}{\sigma^2_n} \bm{v}_n}^2 = \sum_{\sigma_n > \epsilon} |  \ip{\bm{z}}{\bm{v}_n} |^2 \leq \| \bm{z} \|,
}
This gives the second term of \R{coeffics_bounded_indirect}.

For \R{coeffics_conv3}, of course the canonical frame coefficients are well defined since we now assume that $\Phi$ is a frame. We first note that
\bes{
\| \bm{a} - \bm{a}^{\epsilon}_{M,N} \| \leq \sqrt{\sum_{n \in I \backslash I_N} | a_n |^2} + \| \bm{a}_N- \bm{x}^{\epsilon} \|,
}
where $\bm{a}_N = \{ a_n \}_{n \in I_N}$.  Therefore it suffices to consider $\| \bm{a}_N - \bm{x}^{\epsilon} \|$.  Observe that
\bes{
\bm{a}_N = \sum_{n \in I_N} \ip{\bm{a}_N}{\bm{v}_n} \bm{v}_n = \sum_{n \in I_N} \ip{\cS^{-1} f }{\xi_n} \bm{v}_n.
}
Now $\ip{f}{\xi_n}_M=  \ip{\cS \cS^{-1} f}{\xi_n}_M$ and therefore
\eas{
\ip{f}{\xi_n}_M
& = \ip{\cS_N \cS^{-1} f}{\xi_n}_M + \ip{(\cS - \cS_N) \cS^{-1} f}{\xi_n}_M
 = \sigma^2_n \ip{\cT^*_N \cS^{-1} f}{\bm{v}_n} +  \ip{(\cS - \cS_N) \cS^{-1} f}{\xi_n}_M.
}
Notice that $\cS_N \cS^{-1} f \in \rG$ and $(\cS - \cS_N)\cS^{-1} f = f - \cS_N \cS^{-1} f \in \rG$.  Therefore all the terms above are well defined.  Hence, by \R{x_OS_split}, 
\bes{
\bm{x}^{\epsilon} = \sum_{\sigma_{n} > \epsilon} \ip{\cS^{-1} f}{\xi_n} \bm{v}_n + \sum_{\sigma_{n} > \epsilon} \frac{\ip{(\cS - \cS_M) \cS^{-1} f}{\xi_m}_N}{\sigma^2_n} \bm{v}_n ,
}
which gives
\be{
\label{bear1}
\nm{\bm{a}_N-\bm{x}^{\epsilon}} \leq \nm{\sum_{\sigma_n \leq \epsilon} \ip{\cS^{-1} f }{\xi_n} \bm{v}_n} + \nm{ \sum_{\sigma_n > \epsilon} \frac{\ip{(\cS - \cS_M) \cS^{-1} f}{\xi_n}_M}{\sigma^2_n} \bm{v}_n } .
}
Consider the first term.  Let $\bm{z} \in \bbC^N$ be given by $\bm{z} = \sum_{\sigma_n \leq \epsilon}  \ip{\cS^{-1} f }{\xi_n} \bm{v}_n$,
so that the first term is merely $\| \bm{z} \|$.  By the definition of $\lambda^{\epsilon}_{M,N}$, we have $\epsilon \lambda^{\epsilon}_{M,N} \| \bm{z} \| \geq \| \cT_N \bm{z} - \cP^{\epsilon}_{M,N}  \cT_N \bm{z} \|$.  Now, since $\cT_N \bm{z} \perp \rH^{\epsilon}_{M,N}$, we have that $\cP^{\epsilon}_{M,N}  \cT_N \bm{z} = 0$.  Hence
\bes{
\epsilon \lambda^{\epsilon}_{M,N} \| \bm{z} \| \geq \| \cT_N \bm{z} \| = \sup_{\substack{g \in \rH \\ g \neq 0}} \frac{| \ip{\cT_N \bm{z}}{g} |}{\| g \|}.
}
Set $g = \cS^{-1} f$.  Then $\ip{\cT_N \bm{z}}{g} = \sum_{\sigma_n \leq \epsilon} | \ip{\cS^{-1} f}{\xi_m} |^2 = \| \bm{z} \|^2$ and therefore we obtain $\epsilon \lambda^{\epsilon}_{M,N} \| \bm{z} \| \geq \| \bm{z} \|^2/\| \cS^{-1} f \|$.  It follows that
\be{
\label{bear2}
\| \bm{z} \| = \nm{\sum_{\sigma_n \leq \epsilon} \ip{\cS^{-1} f }{\xi_m} \bm{v}_n} \leq \epsilon \lambda^{\epsilon}_{M,N} \nm{\cS^{-1} f } \leq \epsilon \lambda^{\epsilon}_{M,N}/\sqrt{A}  \| \bm{a} \|.
}
Now consider the second term of \R{bear1}.  We have
\eas{
\nm{ \sum_{\sigma_n > \epsilon} \frac{\ip{(\cS - \cS_M) \cS^{-1} f}{\xi_n}_M}{\sigma^2_n} \bm{v}_n }^2  =  \sum_{\sigma_n > \epsilon} \frac{|\ip{(\cS - \cS_N) \cS^{-1} f}{\xi_n}_M |^2}{\sigma^4_n} & \leq \frac{1}{\epsilon^2} \| \cP^{\epsilon}_{M,N} (\cS - \cS_N) \cS^{-1} f \|^2_M
\\
& \leq \frac{1}{\epsilon^2} \|  (\cS - \cS_N) \cS^{-1} f \|^2_M.
}
In the last line, we used Lemma~\ref{l:projection} again. Combining this with \R{bear2} now gives the result.
}

\rem{
\label{r:noisy_data}
These results extend to the setting of noisy measurements. Suppose the measurements are $\bm{y} + \bm{n}$ where $\bm{y} = \cM_M f$ and $\bm{n}$ is a vector of noise. We assume that $\bm{n} \in \mathrm{Ran}(\cM_M)$, in other words it takes the form $\bm{n} = \cM_M g$ for some $g \in \rG$. Then, by linearity, the reconstruction is
\bes{
\tilde{f} = \cP^{\epsilon}_{M,N} f + \cL^{\epsilon}_{M,N} \bm{n},
}
where $\cP^{\epsilon}_{M,N} f$ is the standard reconstruction from the noiseless data $\bm{y}$. Hence, by Theorem \ref{t:indirect_err}, the error satisfies
\be{
\label{noisebound}
\nmu{f - \tilde{f}} \leq \inf_{\bm{z} \in \bbC^N} \left\{  \| f - \cT_N \bm{z} \| + \kappa^{\epsilon}_{M,N} \| f - \cT_N \bm{z} \|_M + \epsilon \lambda^{\epsilon}_{M,N} \| \bm{z} \| \right\} + \kappa^{\epsilon}_{M,N} \nm{\bm{n}}.
}
}
In particular, when $\kappa^{\epsilon}_{M,N}$ is order one, the effect of the noise is proportional to its $\ell^2$-norm.

\subsection{Behaviour of the constants}\label{ss:indirect_constants}

We now consider the behaviour of the constants $\kappa^{\epsilon}_{M,N}$ and $\lambda^{\epsilon}_{M,N}$.  To do so, we define the constant $A'_{M,N}$ as follows:
\be{
\label{discrete_inf}
A'_{M,N} = \inf_{\substack{g \in \rH_N \\ \| g \|=1}} \| g \|^2_M.
}
In general, with a poor choice of sampling functionals, it may be that $A'_{M,N} = 0$. However, even in that case, with assumption \R{data_rich} we have that $\liminf_{M \rightarrow \infty} A'_{M,N} \geq A'$ for any fixed $N$.


\prop{
\label{l:constants_indirect}
The constants $\kappa^{\epsilon}_{M,N}$ and $\lambda^{\epsilon}_{M,N}$ defined in  \R{kappa_lambda_indirect} satisfy
\be{
\label{kappa_lambda_global_indirect}
\kappa^{\epsilon}_{M,N} \leq \frac{\sqrt{B_N}}{\epsilon}, \qquad \lambda^{\epsilon}_{M,N} \leq \frac{\sqrt{B_N}}{\epsilon},
}
for all $M$ and $N$, $M \geq N$, where $B_N$ is the Bessel bound for $\Phi_N$.
Moreover,
\be{
\label{kappa_lambda_limit_indirect}
\kappa^{\epsilon}_{M,N} \leq \frac{1}{\sqrt{A'_{M,N}}},\qquad  \lambda^{\epsilon}_{M,N} \leq \frac{1}{\sqrt{A'_{M,N}}},
}
and if the sampling functionals satisfy \R{data_rich} then, for fixed $N$,
\bes{
\limsup_{M \rightarrow \infty} \kappa^{\epsilon}_{M,N} \leq \frac{1}{\sqrt{A'}} \quad \mbox{and} \quad \limsup_{M \rightarrow \infty} \lambda^{\epsilon}_{M,N} \leq \frac{1}{\sqrt{A'}}.
}
}
\prf{
Let $\bm{y} \in \mathrm{Ran}(\cM_M)$ be given and write $\bm{y} = \cM_M f$ for some $f \in \rG$.  Then, by \R{SVD_proj_def},
\bes{
\| \cL^{\epsilon}_{M,N} \bm{y} \| = \| \cP^{\epsilon}_{M,N} f \| = \nm{\cT_N \sum_{\sigma_n > \epsilon} \frac{\ip{\bm{y}}{\bm{u}_n}}{\sigma_n} \bm{v}_n }.
}
Notice that $\| \cT_N \bm{x} \| \leq \sqrt{B_N} \| \bm{x} \|$. This follows since any frame automatically satisfies the upper Riesz basis condition with constant equal to the Bessel bound.
Hence
\bes{
\| \cP^{\epsilon}_{M,N} f \|^2 \leq B_N \sum_{\sigma_n > \epsilon} \frac{| \ip{\bm{y}}{\bm{u}_n} |^2}{\sigma^2_n} \leq \frac{B_N}{\epsilon^2} \| \bm{y} \|^2,
}
which gives \R{kappa_lambda_global_indirect} for $\kappa^{\epsilon}_{M,N}$.  For \R{kappa_lambda_limit_indirect}, we let $\bm{y} \in \mathrm{Ran}(\cM_M)$ and write $\bm{y} = \cM_M f$ for some $f \in \rG$ once more.
From the definition \R{discrete_inf} of $A'_{M,N}$ and by Lemma~\ref{l:projection} we find
\bes{
\| \cP^{\epsilon}_{M,N} f \| \leq \frac{1}{\sqrt{A'_{M,N}}} \| \cP^{\epsilon}_{M,N} f \|_M \leq \frac{1}{\sqrt{A'_{M,N}}} \nm{f}_{M} =  \frac{1}{\sqrt{A'_{M,N}}} \nm{\bm{y}}.
}
This gives \R{kappa_lambda_limit_indirect}.

We now consider $\lambda^{\epsilon}_{M,N}$.  Let $\bm{z} \in \bbC^M$ be arbitrary.  Using \R{SVD_proj_def} and \R{useful1} we have
\bes{
\cT_N \bm{z} - \cP^{\epsilon}_{M,N}  \cT_N \bm{z} = \cT_N \sum_{\sigma_n \leq \epsilon} \ip{\bm{z}}{\bm{v}_n} \bm{v}_n.  
}
Arguing as above, this implies that $\| \cT_N \bm{z} - \cP^{\epsilon}_{M,N}  \cT_N \bm{z} \|^2 \leq B_N \| \bm{z} \|^2$, which gives \R{kappa_lambda_global_indirect}.  For \R{kappa_lambda_limit_indirect}, we again let $\bm{z} \in \bbC^M$ be arbitrary.  Then
\bes{
\nm{\cT_N \bm{z} - \cP^{\epsilon}_{M,N}  \cT_N \bm{z}}^2_M = \sum_{\sigma_n \leq \epsilon} \sigma^2_n | \ip{\bm{z}}{\bm{v}_n} |^2 \leq \epsilon^2 \| \bm{z} \|^2.
}
Moreover, since $\cT_N \bm{z} - \cP^{\epsilon}_{M,N}  \cT_N \bm{z} \in \rH_N$ we obtain
\bes{
\| \cT_N \bm{z} - \cP^{\epsilon}_{M,N}  \cT_N \bm{z}\|^2 \leq \frac{1}{A'_{M,N}} \| \cT_N \bm{z} - \cP^{\epsilon}_{M,N} \|^2_M \leq \frac{\epsilon^2 }{A'_{M,N}} \| \bm{z} \|^2 ,
}
as required.
}

\subsection{The stable sampling rate}\label{ss:SSR}

Suppose that the sampling functionals satisfy \R{data_rich}. Motivated by Proposition \ref{l:constants_indirect} we now introduce the following concept:

\defn{
\label{d:SSR}
For $1 < \theta < \infty$ and $N \in \bbN$, the \textit{stable sampling rate} is 
\bes{
\Theta^{\epsilon}(N,\theta) = \min \left \{ M \in \bbN : M \geq N,\ \kappa^{\epsilon}_{M,N} \leq \frac{\theta}{\sqrt{A'}},\ \lambda^{\epsilon}_{M,N} \leq \frac{\theta}{\sqrt{A'}} \right \}.
}
}
For a fixed $N$, suppose that $M$ is chosen so that $M \geq \Theta^{\epsilon}(N,\theta)$. Then this guarantees an error bound of the form
\bes{
\| f - \cP^{\epsilon}_{M,N} f \| \leq \inf_{\bm{z} \in \bbC^N} \left\{   \| f - \cT_N \bm{z} \| +\frac{\theta}{\sqrt{A'}} \| f - \cT_N \bm{z} \|_M+ \epsilon \frac{\theta}{\sqrt{A'}} \| \bm{z} \|  \right\} .
}
Hence, sampling according to the stable sampling rate, ensures that the error decays down to roughly $\epsilon$ as $N \rightarrow \infty$. This holds on the additional condition that the term $\| f - \cT_N \bm{z} \|_M \rightarrow 0$; see the discussion after Theorem \ref{t:indirect_err}. Furthermore, sampling according to $ \Theta^{\epsilon}(N,\theta)$ means that the \textit{rate} of decay of the error for finite $N$  depends completely on how well $f$ can be approximated by elements of $\rH_N$ with bounded coefficients. As discussed, this depends completely on the frame $\Phi$ and the element $f$ being approximated. For estimates in certain cases, see \cite{FramesPart1}.

\rem{
If the data is noisy as in Remark \ref{r:noisy_data} and $M \geq \Theta^{\epsilon}(N,\theta)$ then \R{noisebound} becomes
\bes{
\nmu{f - \tilde{f}} \leq \inf_{\bm{z} \in \bbC^N} \left\{   \| f - \cT_N \bm{z} \| +\frac{\theta}{\sqrt{A'}} \| f - \cT_N \bm{z} \|_M+ \epsilon \frac{\theta}{\sqrt{A'}} \| \bm{z} \|  \right\} + \frac{\theta}{\sqrt{A'}} \nm{\bm{n}}.
}  
Note that $\epsilon$ does not enter into the noise term. Recall that the first term will decrease down to a limiting accuracy of at best $\ordu{\epsilon}$. Hence, in the noisy case the limiting accuracy will depend on the maximum of $\epsilon$ and $\nm{\bm{n}}$. In particular, this yields a simple strategy for choosing $\epsilon$ in the noisy case, simply as proportional to the noise level.
}

The behaviour of $\Theta^{\epsilon}(N,\theta)$ as a function of $N$ depends completely on $\Phi$ and the sampling functionals. Thus, theoretical estimates for this quantity can only be established on a case-by-case basis.  We shall consider this issue further in \S \ref{s:SSR} for a particular class of problems. However, there is no general recipe for providing such estimates, and moreover, when possible, doing so typically only reveals the asymptotic growth rate of $\Theta^{\epsilon}(N,\theta)$ with $N$ and not the precise constant.

On the other hand, $\Theta^{\epsilon}(N,\theta)$ can always be computed. To see this, we observe the following:

\lem{
\label{l:kappa_compute}
The constant $\kappa^{\epsilon}_{M,N} $ satisfies
\bes{
\kappa^{\epsilon}_{M,N} \leq \sqrt{\lambda_{\max} \left ( (\bm{B}^{\epsilon}_{M,N})^* \bm{G}_N \bm{B}^{\epsilon}_{M,N} \right ) },
}
where $\bm{B}^{\epsilon}_{M,N} = (\bm{G}^{\epsilon}_{M,N})^{\dag}$ and $\bm{G}_N \in \bbC^{N \times N}$ is the Gram matrix of $\Phi_N$.    If $\mathrm{Ran}(\cM_M) = \bbC^M$ this holds with equality.  The constant $\lambda^{\epsilon}_{M,N}$ satisfies
\bes{
\lambda^{\epsilon}_{M,N} = \epsilon^{-1} \sqrt{\lambda_{\max} \left ( (\bm{C}^{\epsilon}_{M,N})^* \bm{G}_N \bm{C}^{\epsilon}_{M,N} \right )},
}
where $\bm{C}^{\epsilon}_{M,N} = \bm{V} (\bm{I} - \bm{I}^{\epsilon}) \bm{V}^*$ and $\bm{I}_{\epsilon}$ is the diagonal matrix with $n^{\rth}$ entry $1$ if $\sigma_n \geq \epsilon$ and zero otherwise.
}
\prf{
Let $\bm{y} \in \mathrm{Ran}(\cM_M)$ with $\nm{\bm{y}} = 1$.  Then
\bes{
\nm{\cL^{\epsilon}_{M,N} \bm{y}}^2 = \nm{\cT_N \bm{B}^{\epsilon}_{M,N} \bm{y} }^2 = \bm{y}^* (\bm{B}^{\epsilon}_{M,N})^* \bm{G}_N \bm{B}^{\epsilon}_{M,N} \bm{y} \leq \lambda_{\max} \left ( (\bm{B}^{\epsilon}_{M,N})^* \bm{G}_N \bm{B}^{\epsilon}_{M,N} \right ),
}
since $\bm{G}_N = \cT^*_N \cT_N$.  This gives the first result.

Let $\bm{z} \in \bbC^N$, $\nm{\bm{z}}=1$.  Then
\eas{
\nm{\cT_N \bm{z} - \cP^{\epsilon}_{M,N} \cT_N \bm{z}}^2 = \nm{\cT_N \left ( \bm{I}- (\bm{G}^{\epsilon}_{M,N})^{\dag} \bm{G}_{M,N} \right ) \bm{z} }^2  &= \nm{\cT_N \left ( \bm{I} - \bm{V} \bm{I}^{\epsilon} \bm{V}^* \right ) \bm{z}}^2 
\\
&= \nm{\cT_N \bm{C}^{\epsilon}_{M,N} \bm{z} }^2 = \bm{z}^* (\bm{C}^{\epsilon}_{M,N})^* \bm{G}_N \bm{C}^{\epsilon}_{M,N} \bm{z}.
}
Maximizing over $\bm{z}$ now gives the result.
}

This result implies that $\kappa^{\epsilon}_{M,N}$ and $\lambda^{\epsilon}_{M,N}$ can be computed, and therefore so can $\Theta^{\epsilon}(N;\theta)$, provided the matrix $\bm{G}_N$ has been computed.

\rem{
\label{r:SSRcompute}
In practice, it may be difficult to compute $\bm{G}_N$, since its entries are inner products which may for instance be integrals.  This may be overcome by a further approximation, e.g.\ a quadrature.  Specifically, let $K \geq 1$ and $\jmath_{k,K}$ be a family of linear functionals such that
\bes{
\lim_{K \rightarrow \infty} \sum^{K}_{k=1} \jmath_{k,K}(f) \overline{\jmath_{k,K}(g)} = \ip{f}{g},\qquad \forall f,g \in \rG.
}
Let $\bm{H}_{K,N} = \{ \jmath_{k,K}(\phi_n) \}^{K}_{k=1, n \in I_N} \in \bbC^{K \times N}$.  Then $(\bm{H}_{K,N})^* \bm{H}_{K,N} \approx \bm{G}_N$ for large $K$.  Hence, by the previous lemma (assuming $\mathrm{Ran}(\cM_M) = \bbC^M$ for ease of presentation), we have
\bes{
\kappa^{\epsilon}_{M,N} \approx \nm{\bm{H}_{K,N} \bm{B}^{\epsilon}_{M,N}}_{2} = \nm{\bm{H}_{K,N} \bm{V} (\bm{\Sigma}^{\epsilon})^{\dag} }_{2},\qquad  \lambda^{\epsilon}_{M,N} \approx \nm{\bm{H}_{K,N} \bm{V} ( \bm{I} - \bm{I}^{\epsilon} ) }_{2},
}
for sufficiently large $K$.  If, for instance, the functionals $\jmath_{k,K}$ correspond to pointwise evaluations as part of a quadrature, this gives a means of numerically approximating $\kappa^{\epsilon}_{M,N}$ and $\lambda^{\epsilon}_{M,N}$.
}

\section{Frame approximation from frame samples}\label{s:framesforall}

In this section, we discuss Example \ref{ex:otherframe}, in which both the approximation system and sampling functionals arise from frames of $\rH$. We write $\Phi = \{ \phi_n \}_{n \in I}$ for the approximation frame (with bounds $A$ and $B$) and $\Psi = \{ \psi_m \}_{m \in J}$ for the sampling frame (with bounds $A'$ and $B'$). We assume that $\{ I_N \}_{N \in \bbN}$ is a sequence of nested index sets with
\bes{
I_N \subset I,\ |I_N| = N,\ \forall N \in \bbN,\qquad \bigcup^{\infty}_{N=1} I_N = I,
}
and similarly, we assume that $\{ J_M \}_{M \in \bbN}$ is a sequence of nested index sets with
\bes{
J_M \subset I,\ |J_M| = M,\ \forall M \in \bbN,\qquad \bigcup^{\infty}_{M=1} I_M = J.
}
As before, we write $\cP^{\epsilon}_{M,N} f$ for the truncated SVD approximation of $f \in \rH$ (note that $\rG = \rH$ in this case, since the sampling functionals arise from a frame of $\rH$). We write $\bm{x}^{\epsilon}$ for its coefficients.

\subsection{Error and coefficient bounds}

\thm{
In the above setting, the truncated SVD projection $\cP^{\epsilon}_{M,N} f$ of $f \in \rH$ satisfies
\be{
\label{err_bd_frame_samp}
\| f - \cP^{\epsilon}_{M,N}f \| \leq \inf_{\bm{z} \in \bbC^N} \left\{ \left ( 1 + \sqrt{B'} \kappa^{\epsilon}_{M,N} \right ) \nm{f - \cT_{N} \bm{z}} + \epsilon \lambda^{\epsilon}_{M,N} \nm{\bm{z}} \right\}.
}
Its coefficients $\bm{x}^{\epsilon}$ satisfy
\be{
\label{coeffics_bounded2}
\| \bm{x}^{\epsilon} \| \leq \inf_{\bm{z} \in \bbC^N} \left\{ \frac{\sqrt{B'}}{\epsilon} \nm{f - \cT_{N} \bm{z} } + \nm{\bm{z}} \right\},
}
and, if $\bm{a}^{\epsilon}_{M,N} \in \ell^2(I)$ the extension of $\bm{x}^\epsilon$ by zero,
\be{
\label{coeffics_conv2}
\| \bm{a} - \bm{a}^{\epsilon}_{M,N} \| \leq \left ( 1 + \frac{\sqrt{B B'}}{\epsilon} \right ) \sqrt{\sum_{n \in I \backslash I_N} | a_n |^2} + \epsilon \frac{\lambda^{\epsilon}_{M,N}}{\sqrt{A}} \| \bm{a} \|,
}
where $\bm{a} =  \{ \ip{f}{\cS^{-1}\phi_n} \}_{n \in I}$ are the canonical frame coefficients of $f \in \rH$.
}
\prf{
Recall \R{framesamplesUB}. The first two bounds now follow immediately from Theorems \ref{t:indirect_err} and \ref{t:coefficients_indirect} respectively. For the third bound, we use use \R{coeffics_conv3} and then observe that
\bes{
\nm{ (\cS - \cS_N) \cS^{-1} f }^2 = \nm{\sum_{n \in I \backslash I_N} a_n \phi_n }^2 \leq B \sum_{n \in I \backslash I_N} |a_n |^2,
}
where for the final step we recall that a frame with upper frame bound $B$ satisfies the upper Riesz basis condition with constant $B$. The result now follows immediately.
}

\subsection{Behaviour of the coefficients and $\ord{\epsilon}$ accuracy}

We now consider the constants $\kappa^{\epsilon}_{M,N}$ and $\lambda^{\epsilon}_{M,N}$ and the stable sampling rate:

\prop{
\label{l:constants}
The constants $\kappa^{\epsilon}_{M,N}$ and $\lambda^{\epsilon}_{M,N}$ satisfy
\be{
\label{kappa_lambda_global}
\kappa^{\epsilon}_{M,N} , \lambda^{\epsilon}_{M,N} \leq \left \{ \begin{array}{cc} \sqrt{B}/\epsilon & \Psi \neq \Phi \\ 1/\sqrt{\epsilon} & \Psi = \Phi,
\end{array} \right . ,
}
for all $M$ and $N$, $M \geq N$.
Moreover, for fixed $N$,
\bes{
\limsup_{M \rightarrow \infty} \kappa^{\epsilon}_{M,N} \leq \frac{1}{\sqrt{A'}} \quad \mbox{and} \quad \limsup_{M \rightarrow \infty} \lambda^{\epsilon}_{M,N} \leq \frac{1}{\sqrt{A'}},
}
}

This result is essentially a special case of Proposition \ref{l:constants_indirect}, except in the case where $\Psi = \Phi$ where we have a slightly improved worst-case behaviour, with the right-hand side of \R{kappa_lambda_global} scaling like $1/\sqrt{\epsilon}$ as opposed to $1/\epsilon$. As is made clear by the proofs, this discrepancy is due to the fact that in the latter case the measurements are just inner products with respect to the same frame. We prove this in a moment. First, however, we consider its implications for limiting accuracy:

\cor{
\label{c:frame_oversampling}
For each $1 < \theta < \infty$ the truncated SVD approximation $\cP^{\epsilon}_{M,N} f$ satisfies
\bes{
\limsup_{\substack{M,N \rightarrow \infty \\ M \geq \Theta^{\epsilon}(N,\theta)}}\| f - \cP^{\epsilon}_{M,N} f \| \leq \epsilon \frac{\theta}{\sqrt{A A'}} \| f \|,
}
where $\Theta^{\epsilon}(N,\theta)$ is as in Definition \ref{d:SSR}.
Moreover, the coefficients $\bm{x}^{\epsilon}$ satisfy
\bes{
\limsup_{\substack{M,N \rightarrow \infty \\ M \geq \Theta^{\epsilon}(N,\theta)}} \nm{\bm{x}^{\epsilon}} \leq \frac{1}{\sqrt{A}} \nm{f},\qquad \limsup_{\substack{M,N \rightarrow \infty \\ M \geq \Theta^{\epsilon}(N,\theta)}} \nm{\bm{a} - \bm{a}^{\epsilon}_{M,N}} \leq \epsilon \frac{\theta}{\sqrt{AA'}} \nm{\bm{a}}.
}
}
\prf{
The proof is based on the canonical frame coefficients $\bm{a} = \{ \ip{f}{\cS^{-1} \phi_n } \}_{n \in I}$.  Let $\bm{z} = \{ a_n \}_{n \in I_N}$.  Then $\nm{\bm{z}} \leq \nm{\bm{a}} \leq 1/\sqrt{A} \nm{f}$ since the dual frame has upper frame bound $A^{-1}$ (see \S \ref{ss:canonical_dual_frame}).  Therefore \R{err_bd_frame_samp} and Proposition \ref{l:constants} gives
\bes{
\| f - \cP^{\epsilon}_{M,N} f \| \leq \left(1+\sqrt{\frac{B}{A}} \theta \right)  \nm{f - \sum_{n \in I_N} \ip{f}{\cS^{-1} \phi_n} \phi_n } + \epsilon \frac{\theta}{A}\nm{f} .
}
As $N \rightarrow \infty$ \R{dual_rep} gives that the first term vanishes.  Hence we obtain the result for $f$.  For the other results, we use \R{coeffics_bounded2} and \R{coeffics_conv2} instead.
}

In summary, provided $M$ is chosen above the stable sampling rate $\Theta^{\epsilon}(N,\theta)$, the approximation $\cP^{\epsilon}_{M,N} f$ converges to within roughly $\epsilon$ of $f$ and the coefficients converge to within roughly $\epsilon$ of the frame coefficients $\bm{a}$, and in particular are small in norm for large $N$. As a consequence, in the setting of Theorem \ref{t:Proj_err} where $\Psi = \Phi$, we overcome the $\sqrt{\epsilon}$ bottleneck, with the limiting accuracy bounded by $\epsilon \theta \nm{f} / A$ as opposed to $\sqrt{\epsilon} \nm{f} / \sqrt{A}$ (see Corollary \ref{c:frame_projection}).

This result also illuminates the role that the frame structure plays in both the approximation and the sampling. Indeed, the limiting error depends on both of the lower frame bounds $A$ and $A'$, while the limiting size of the coefficients depends only on $A$. This is as expected. The existence of small norm coefficients depends only on the approximation frame $\Phi$, a small limiting error depends on both the sampling frame $\Psi$ and the approximation frame $\Phi$.

\prf{
[Proof of Proposition \ref{l:constants}]

All results follow immediately from Proposition \ref{l:constants_indirect}, except for \R{kappa_lambda_global} in the case $\Psi = \Phi$ for which we require a different argument.

Consider $\kappa^{\epsilon}_{M,N}$ first.  Let $\bm{y} \in  \mathrm{Ran}(\cM_M)$ be given and notice that we may write $\bm{y} = \cT^*_M f$ for some $f \in \rH$ so that $
\| \cL^{\epsilon}_{M,N} \bm{y} \| = \| \cP^{\epsilon}_{M,N} f \|$.  By \R{x_SVD_def} we have
\bes{
\| \cP^{\epsilon}_{M,N} f \|^2 = \sum_{\sigma_m,\sigma_n > \epsilon}  \frac{\ip{\bm{y}}{\bm{u}_m} \overline{\ip{\bm{y}}{\bm{u}_n} } }{\sigma_m \sigma_n} \ip{\xi_m}{\xi_n}.
}
Recall that $\cM_M \cT_N \bm{v}_m = \sigma_m \bm{u}_m$. Since $\cM_M = \cT^*_M$ in this case, we have $\cT^*_N \cT_N \bm{v}_m = \sigma_m \bm{\tilde{u}}_m$, where where $\bm{\tilde{u}}_m \in \bbC^N$ is the vector with entries $(\bm{\tilde{u}}_m)_k = (\bm{u}_m)_k$ for $k \in I_N$. Hence $\ip{\xi_m}{\xi_n} = \ip{\cT^*_N \cT_N \bm{v}_m}{\bm{v}_n} = \sigma_m \ip{\bm{\tilde{u}}_m}{\bm{v}_n}$ and this gives
\eas{
\| \cP^{\epsilon}_{M,N} f \|^2 = \ipl{ \sum_{\sigma_n > \epsilon} \ip{\bm{y}}{\bm{u}_m} \bm{\tilde{u}}_m  }{  \sum_{\sigma_n > \epsilon}  \frac{\ip{\bm{y}}{\bm{u}_n}}{\sigma_n} \bm{v}_n  }
\leq \nm{\sum_{\sigma_n > \epsilon} \ip{\bm{y}}{\bm{u}_m} \bm{\tilde{u}}_m } \nm{ \sum_{\sigma_n > \epsilon}  \frac{\ip{\bm{y}}{\bm{u}_n}}{\sigma_n} \bm{v}_n  }.
}
By orthogonality, the second term satisfies
\bes{
\nm{ \sum_{\sigma_n > \epsilon}  \frac{\ip{\bm{y}}{\bm{u}_n}}{\sigma_n} \bm{v}_n  } =  \sqrt{\sum_{\sigma_n > \epsilon} \frac{| \ip{\bm{y}}{\bm{u}_n} |^2}{\sigma^2_n}} \leq \frac{\nm{\bm{y}}}{\epsilon}.
}
Consider the first term.  Let $\cQ_N : \bbC^M \rightarrow \bbC^M$ be the projection defined by $(\cQ_N \bm{x})_{m} = x_m$, $m \in I_N$ and $(\cQ_N \bm{x})_m = 0$, $m \in I_M \backslash I_N$.  Then
\bes{
\nm{\sum_{\sigma_n > \epsilon} \ip{\bm{y}}{\bm{u}_m} \bm{\tilde{u}}_m } = \nm{\cQ_N \left (\sum_{\sigma_n > \epsilon} \ip{\bm{y}}{\bm{u}_m} \bm{u}_m  \right ) } \leq \nm{\sum_{\sigma_n > \epsilon} \ip{\bm{y}}{\bm{u}_m} \bm{u}_m} \leq \nm{\bm{y}}.
}
Therefore, we deduce that 
\bes{
\| \cP^{\epsilon}_{M,N} f \|^2 \leq \nm{\bm{y}}/\epsilon,
}
and the result for $\kappa^{\epsilon}_{M,N}$ now follows from its definition.

Now consider $\lambda^{\epsilon}_{M,N}$.  Let $\bm{z} \in \bbC^N$ be arbitrary and recall that
\bes{
\cT_N \bm{z} - \cP^{\epsilon}_{M,N}  \cT_N \bm{z} = \sum_{\sigma_n \leq \epsilon} \ip{\bm{z}}{\bm{v}_n} \xi_n.
}
Hence
\bes{
\nm{\cT_N \bm{z} - \cP^{\epsilon}_{M,N}  \cT_N \bm{z}}^2 = \sum_{\sigma_m,\sigma_n \leq \epsilon} \ip{\bm{z}}{\bm{v}_n} \overline{\ip{\bm{z}}{\bm{v}_n}} \ip{\xi_m}{\xi_n}.
}
As above, we note that $\ip{\xi_m}{\xi_n} = \sigma_m \ip{\bm{\tilde{u}}_m}{\bm{v}_n}$, and therefore
\bes{
\nm{\cT_N \bm{z} - \cP^{\epsilon}_{M,N}  \cT_N \bm{z}}^2  = \ipl{ \sum_{\sigma_n \leq \epsilon} \sigma_m \ip{\bm{z}}{\bm{v}_m} \bm{\tilde{u}}_m }{\sum_{\sigma_n \leq \epsilon} \ip{\bm{z}}{\bm{v}_n} \bm{v}_n} \leq \epsilon \| \bm{z} \|^2.
}
Since $\bm{z}$ was arbitrary, we now obtain the result for $\lambda^{\epsilon}_{M,N}$.
}

\section{ONB$+1$ and ONB$+K$ frames}\label{s:SSR}




We conclude this paper with several examples to illustrate the stable sampling rate.   Let $\{ \varphi_{n} \}_{n \in \bbN}$ be an orthonormal basis of $\rH$ and $\psi \in \rH$, $\nm{\psi} = 1$, be such that $\ip{\psi}{\varphi_n} \neq 0$ for infinitely-many $n \in \bbN$.  Then the indexed family
\bes{
\Phi = \{ \phi_0,\phi_1,\ldots \} = \{ \psi , \varphi_1 , \varphi_2 ,\ldots \},
}
is a frame for $\Phi$ with frame bounds $A = 1$ and $B = 2$. We refer to this frame as the \textit{ONB$+1$ frame}.  Note that it was previously used in \cite{FramesPart1} to show that the Gram matrix of a frame can be arbitrarily badly conditioned.  It is motivated by the idea of 
`enriching' an orthonormal basis to better capture a certain feature of a function under approximation (e.g.\ a singularity or oscillation).

Throughout this section, we let $\cQ_{N}$ denote the projection onto $\spn \{ \varphi_1,\ldots,\varphi_{N} \}$, i.e.
\bes{
\cQ_{N} f = \sum^{N}_{n=1} \ip{f}{\varphi_n} \varphi_n.
}

\subsection{The stable sampling rate for the ONB$+1$ frame}

A problem of interest is that where the samples are inner products with respect to the orthonormal basis $\{ \varphi_m \}_{m \in \bbN}$.  That is,
\be{
\label{sampleONB}
\ell_{m,M}(f) = \ell_m(f) = \ip{f}{\varphi_m},\quad m = 1,\ldots,M,\ M \in \bbN.
}
For instance, these are Fourier coefficients if $\{ \varphi_m \}_{m \in \bbN}$ is the Fourier basis, and hence the goal would be to compute a better approximation in the frame $\Phi$ from the given Fourier data.
Note that this is an instance of \S \ref{s:framesforall} with $A' = B' = 1$.  Note also that $\nm{g}_M = \nm{\cQ_M g}$.  Recalling Proposition \ref{l:constants_indirect}, we determine the stable sampling rate we note that it suffices to estimate
\be{
\label{infsupbd}
A'_{M,N}= \inf_{\substack{g \in \rH_N \\ \nm{g} = 1}} \nm{\cQ_M g}^2 = 1 - \sup_{\substack{g \in \rH_N \\ \nm{g} = 1}} \nm{g - \cQ_M g}^2,
}
where $\rH_N = \spn \{ \psi , \varphi_1,\ldots,\varphi_{N-1} \}$.

\lem{
\label{l:ONBplus1indirect}
For $M \geq N$, we have
\bes{
\sup_{\substack{g \in \rH_N \\ \nm{g} = 1}} \nm{g - \cQ_M g} = \frac{\nm{\psi - \cQ_M \psi}}{\nm{\psi - \cQ_{N-1} \psi}}.
}
}
\prf{
Let $g \in \rH_N$ and write $g = x_0 \psi + \sum^{N-1}_{n=1} x_{n} \varphi_{n}$.  Then
\bes{
\ip{g}{\varphi_n} = x_0 \ip{\psi}{\varphi_n} + x_{n},\quad n = 1,\ldots,N-1.
}
Therefore
\bes{
g = x_0 \psi + \sum^{N-1}_{n=1} \left ( \ip{g}{\varphi_n} - x_0 \ip{\psi}{\varphi_n} \right ) \varphi_n = x_0 (\psi - \cQ_{N-1} \psi )+ \cQ_{N-1} g
}
Rearranging gives
\bes{
g - \cQ_{N-1} g = x_0 (\psi - \cQ_{N-1} \psi ),
}
and taking the norm of both sides, we find that
\bes{
|x_0| = \frac{\nm{g-\cQ_{N-1} g}}{\nm{\psi - \cQ_{N-1} \psi}} \leq \frac{\nm{g}}{\nm{\psi - \cQ_{N-1} \psi}} = \frac{1}{\nm{\psi - \cQ_{N-1} \psi}}.
}
Also, we have
\bes{
\ip{g}{\varphi_n} = x_0 \ip{\psi}{\varphi_n},\quad n \geq N.
}
Therefore
\bes{
\nm{g - \cQ_{M} g} = |x_0| \nm{\psi - \cQ_M \psi}.
}
Combining this with the bound for $|x_0|$ gives that
\bes{
\sup_{\substack{g \in \rH_N \\ \nm{g} = 1}} \nm{g - \cQ_M g} \leq \frac{\nm{\psi - \cQ_M \psi}}{\nm{\psi - \cQ_{N-1} \psi}}
}
To show equality, we divide into two cases. Suppose first that $\psi \perp \spn \{ \varphi_1,\ldots,\varphi_{N-1}\}$. Then $\cQ_{N-1} \psi = 0$ and so we may take $g = \psi / \nm{\psi}$ to obtain equality. On the other hand, if $\psi \not\perp \spn \{ \varphi_1,\ldots,\varphi_{N-1}\}$ then there exists a $g \in \rH_N$, $\nm{g} = 1$, with $\cQ_{N-1} g = 0$. In this case, the above arguments give that $\nm{g - \cQ_{M} g} = \nm{\psi - \cQ_M \psi} / \nm{\psi - \cQ_{N-1} \psi}$, which implies the result.
}

This leads us to the following result:

\thm{
\label{t:ONBplus1SSR}
Suppose that $\psi$ is such that $| \ip{\psi}{\varphi_n} | \sim c n^{-\alpha}$ as $n \rightarrow \infty$ for some $c>0$ and $\alpha > 1/2$. Then the stable sampling rate
\bes{
\Theta^{\epsilon}(N,\theta) \leq C N,
}
for some constant $C>0$ depending on $c$, $\alpha$ and $\theta$ only.  Conversely, if $| \ip{\psi}{\varphi_n} | \sim c \rho^{-n}$ as $n \rightarrow \infty$ for some $c>0$ and $\rho > 1$ then
\bes{
\Theta^{\epsilon}(N,\theta) \leq N + C,
}
where $C>0$ depends on $c$, $\rho$ and $\theta$ only.
}
\prf{
In the first case, the condition on the coefficients gives
\bes{
\nm{\psi - \cQ_{N} \psi}^2 = \sum_{n \geq N} | \ip{\psi}{\varphi_n} |^2 \sim c' N^{1-2\alpha},\quad N \rightarrow \infty,
}
where $c'$ depends on $c$ and $\alpha$.
Hence Lemma \ref{l:ONBplus1indirect} and the bound \R{infsupbd} give
\bes{
A'_{M,N} \geq 1 - c'' \left ( \frac{M^{1/2-\alpha} }{N^{1/2-\alpha}} \right )^2,
}
for some constant $c''$ depending on $c$ and $\alpha$. Recalling that $B' = A' = A = 1$ and $B = 2$ for this frame and using Proposition \ref{l:constants_indirect} gives the first result. For the second result, we notice that 
\bes{
\nm{\psi - \cQ_{N} \psi}^2 \sim c \frac{\rho^{-N}}{1-\rho}.
}
We now argue as in the previous case.
}


This result shows that the stable sampling rate is linear when $\psi$ has algebraically or exponentially-decaying coefficients in the orthonormal basis $\{ \varphi_n \}_{n \in \bbN}$.  Furthermore, the better $\psi$ is approximated in this basis, the smaller the stable sampling rate is, as evidenced by the case of exponentially-decaying coefficients.
In fact, Lemma \ref{l:ONBplus1indirect} demonstrates the connection between the stable sampling rate and how well approximated $\psi$ is in the orthonormal basis $\{ \varphi_n \}_{n \in \bbN}$.  Specifically, the faster the projection errors $\nm{\psi - \cQ_M \psi}$ decay, the smaller $M \geq N$ needs to be so that $\frac{\nm{\psi - \cQ_M \psi}}{\nm{\psi - \cQ_{N-1} \psi}} \leq \delta$ for constant $0 < \delta < 1$.  This is intuitive.  The better $\psi$ is approximated by this basis, the more information the data, i.e.\ inner products with the $\varphi_n$, carries about the element $g$.

On the other hand, the worse $g$ is approximated the higher the stable sampling rate.  Indeed, if $\nm{\psi - \cQ_M \psi} \asymp (\log(N))^{-1}$ then it is a simple exercise to show that the stable sampling rate is algebraic in $N$ with the power depending on $\theta$, i.e.\ $\Theta^{\epsilon}(N,\theta) = \ord{N^{h(\theta)}}$ for some function $h(\theta) \geq 1$ with $h(\theta) \rightarrow \infty$ as $\theta \rightarrow 1^+$.  


\rem{
One can also determine a bound on the stable sampling rate for special case $\Phi = \Psi$, also discussed in \S \ref{s:framesforall}. In this case, the data is given by the inner products
\bes{
\ip{f}{\psi},\quad \ip{f}{\varphi_m},\quad m=1,\ldots,M-1.
}
Indeed, observe that
\bes{
A'_{M,N} = \inf_{\substack{g \in \rH_N \\ \nm{g} = 1}} \left ( | \ip{g}{\psi} |^2 + \nm{\cQ_{M-1} g}^2 \right ) \geq  \inf_{\substack{g \in \rH_N \\ \nm{g} = 1}}  \nm{\cQ_{M-1} g}^2.
}
The right-hand side is precisely the constant in \R{infsupbd} with $M$ replaced by $M-1$.  Hence, up to an additive factor of one, the stable sampling rate for this problem satisfies the same bounds as those of Theorem \ref{t:ONBplus1SSR}.
}

\subsection{The approximation of functions with logarithmic singularities}

Let $\rH = L^2(0,1)$.  The scaled Legendre polynomials, $\varphi_{n}(x) = \sqrt{2n-1} P_{n-1}(2x-1)$, $n \in \bbN$, form an orthonormal basis for $\rH$.  Here $P_n(x)$ is the usual Legendre polynomial, with normalization $P_n(1) = 1$.  This basis is extremely good at approximating smooth functions.  However, many functions that may arise in applications, such as Green's functions or solutions to PDEs on domains with corners, fail to be smooth at a point $x$, yet posses a known type of singularity there. That is, in these applications  we may want to approximate functions of the form
\be{
\label{faug}
f(x) = g(x) + w(x) h(x),\quad x \in (0,1)
}
where $g,h$ are smooth functions, and $w \in L^2(0,1)$ is a known function which may be singular at, say, $x = 0$.  Such functions cannot generally be accurately approximated using polynomials alone.  However, they can be more accurately captured by enriching the polynomial basis with the function $w$.  This gives a frame
\be{
\label{Phiaug1}
\Phi = \{ \varphi_{n} \}^{\infty}_{n=1} \cup \{ w \},
}
for $\rH$.  Indeed, since the $\varphi_n$ are an orthonormal basis, it quickly follows that
\bes{
\nm{f}^2_{2} = \sum^{\infty}_{n=1} | \ip{f}{\varphi_{n}} |^2 \leq \sum^{\infty}_{n=1} | \ip{f}{\varphi_{n}} |^2 + | \ip{f}{w} |^2 \leq \nm{f}^2 + \nm{f}^2 \nm{w}^2 .
}
Hence this is a frame with bounds $A \geq 1$ and $B \leq 1 + \nm{w}^2$.

\setlength{\figurewidth}{7.5cm}
\setlength{\figureheight}{5.5cm}

\begin{figure}
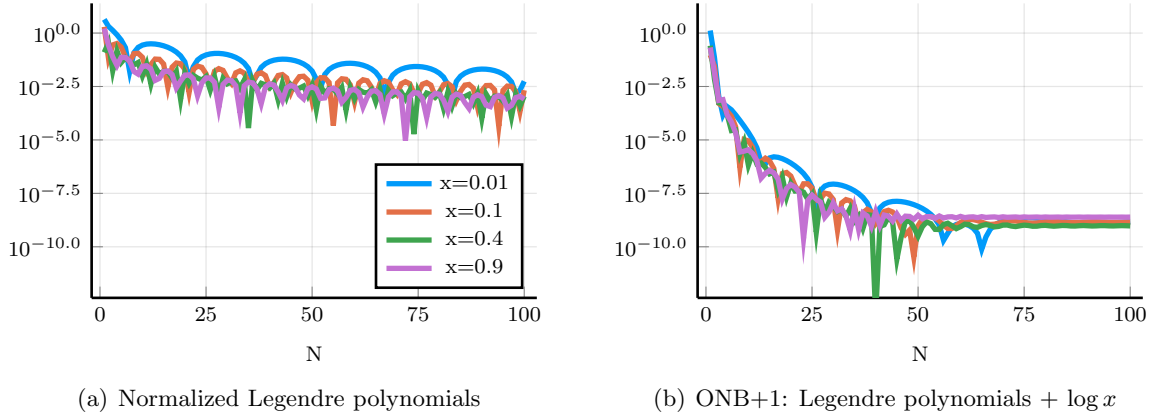

\begin{center}
\subfigure[Normalized Legendre polynomials]{
\begin{small}
\input{legendre_continuous_rev.tikz}
\end{small}
}
\subfigure[ONB$+1$: Legendre polynomials + $\log x$]{
\begin{small}
\input{legendre_log_continuous_rev.tikz}
\end{small}
}
\caption{Pointwise error as a function of the polynomial degree $N$ for the approximation of the logarithmically singular function \eqref{f_example} on $[0,1]$ using Legendre polynomials (left panel) and Legendre polynomials augmented with $\log x$ (right panel).  The error is shown in four points in the interval $[0,1]$. In both cases, the approximation problem is solved using generalized sampling \eqref{sampleONB} with $M=2N$. The generalized samples $(\langle f,\varphi_m\rangle )_{m=1}^{M}$ were evaluated using adaptive numerical integration. The regularization threshold is $\epsilon = 2e^{-13}$.}\label{fig:figureA}
\end{center}
\end{figure}

The case of a logarithmic singularity, i.e.\ $w(x) = \log(x)$, is an important instance of the problem.  Fig.\ \ref{fig:figureA} gives an illustration of the benefits of this frame over just the polynomial basis for approximating the simple yet singular function
\be{
\label{f_example}
f(x) = e^x + \log (x) \cos (x).
}
The polynomial interpolation to $f$ converges poorly, as expected. However, adding just the single element $w(x)=\log(x)$ to the basis results in significantly faster convergence rates, shown in Fig.\ \ref{fig:figureA}(b). Importantly, note that the approximation scheme does not evaluate the smooth parts of $f$ separately. They are implicitly approximated simultaneously when approximating $f$ from its samples. Indeed, if the smooth parts of $f$ were known separately in an application, the approximation problem simplifies and there would be no need to construct a frame. Note also that the evaluation of the generalized samples \eqref{sampleONB} requires the evaluation of integrals, and this step is computationally demanding because the integrals are weakly singular. In subsequent examples, we shall consider a fully discrete approximation based on function samples.

Using Theorem \ref{t:ONBplus1SSR}, we may estimate the stable sampling rate for this problem:

\prop{
Let $\rH = L^2(0,1)$, $w(x) = \log(x)$, $\{ \varphi_n \}$ be the Legendre basis on $\rH$, $\Phi$ be as in \R{Phiaug1} and consider the sampling functionals \R{sampleONB}.  Then the stable sampling rate for this problem is linear in $N$, and specifically,
\bes{
\Theta^{\epsilon}(N,\theta) \leq \max \left \{ N , \frac{N-1}{\sqrt{1-1/\theta^2}} \right \},\quad \forall \theta > 1,\ N \geq 2.
}
}
\prf{
The Legendre polynomials satisfy
\bes{
\int^{1}_{0} P_n(2x-1) \log(x) \D x = \frac{(-1)^{n+1}}{n(n+1)},\quad n \geq 1.
}
This follows from the differential equation $((1-x^2) P'_n(x))' + n(n+1) P_n(x) = 0$ after two integrations by parts.
Let $\psi(x) = \log(x)$.  Then, for $M \geq 1$,
\bes{
\nm{\psi - \cQ_{M} \psi}^2 = \sum_{m > M} | \ip{\psi}{\varphi_m} |^2 = \sum_{m \geq M} \frac{2m+1}{m^2(m+1)^2}  = \sum_{m \geq M} \left ( \frac{1}{m^2} - \frac{1}{(m+1)^2} \right ) = \frac{1}{M^2}.
}
Lemma \ref{l:ONBplus1indirect} now gives
\bes{
 \sup_{\substack{g \in \rH_N \\ \nm{g} = 1}} \nm{g - \cQ_M g} = \frac{N-1}{M}.
}
Therefore 
\bes{
 \sup_{\substack{g \in \rH_N \\ \nm{g} = 1}} \nm{g - \cQ_M g} \leq \sqrt{1-1/\theta^2},
}
provided 
\bes{
\frac{N-1}{M} \leq \sqrt{1-1/\theta^2} \qquad \Leftrightarrow \qquad
M \geq \frac{N-1}{\sqrt{1-1/\theta^2}},
}
as required.
}

\subsection{ONB+$K$ frames}

Functions with logarithmic singularities can be more acccurately approximately using the frame \R{Phiaug1} than the Legendre polynomial basis alone.  However, the accuracy may be limited, due to the presence of weak logarithmic singularities.  To increase the accuracy one may consider a frame of the type
\be{
\label{Phiaug}
\Phi = \{ \varphi_{n} \}^{\infty}_{n=1} \cup \{ \psi_k \}^{K}_{k=1},
}
for fixed $K \geq 1$,  where $\psi_k(x) = w(x) \varphi_k(x)$.  

\prop{
Let $\{ \varphi_n \}$ be the Legendre basis on $\rH$ and $w \in L^2(0,1)$.  Then \R{Phiaug} is a frame for any fixed $K \geq 1$, with frame bounds
\bes{
1 \leq A \leq B \leq 1 + \nm{w}^2 K^2.
}
}
\prf{
First observe that $\psi_k \in L^2(0,1)$ since $w \in L^2(0,1)$ and $\varphi_{k} \in L^{\infty}(0,1)$.  Second, we have
\bes{
\nm{f}^2 = \sum^{\infty}_{n=1} | \ip{f}{\varphi_n} |^2 \leq \sum^{\infty}_{n=1} | \ip{f}{\varphi_n} |^2 + \sum^{K}_{k=1} | \ip{f}{\psi_k} |^2 ,
}
and therefore the lower frame condition holds with $A \geq 1$.  Moreover,
\eas{
\sum^{\infty}_{n=1} | \ip{f}{\varphi_n} |^2 + \sum^{K}_{k=1} | \ip{f}{\psi_k} |^2  &\leq \nm{f}^2 + \nm{f}^2 \nm{w}^2 \sum^{K}_{k=1} \nm{\varphi_k}^2_{L^{\infty}} 
\\
& = \nm{f}^2 \left ( 1 + \nm{w}^2 \sum^{K}_{k=1} (2k-1) \right ) = \nm{f}^2 \left ( 1 + \nm{w}^2 K^2 \right ).
}
Here, in the penultimate step, we use the fact that $| \varphi_{n}(x) | \leq \varphi_n(1) = \sqrt{2n-1}$ for $0 \leq x \leq 1$.  This completes the proof.
}

\begin{figure}[t]
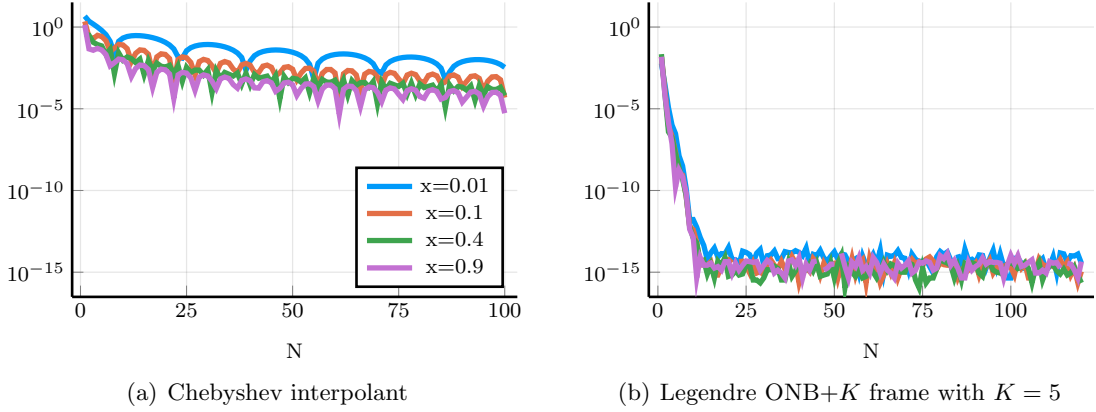

\begin{center}
\subfigure[Legendre interpolation]{
\begin{tiny}
\input{legendre_interpolation_rev.tikz}
\end{tiny}
}
\subfigure[Legendre ONB$+K$ frame with $K=5$]{
\begin{tiny}
\input{legendre_log5_discrete_legendre_rev.tikz}
\end{tiny}
}
\caption{Pointwise error as a function of the polynomial degree $N$ for the approximation of the logarithmically singular function \eqref{f_example} on $[0,1]$, using Legendre interpolation (left panel) and the ONB$+K$ frame $\Phi$ with Legendre polynomials, $K=5$  and $w(x)=\log x$ (right panel). The error is shown in three points in the interval $[0,1]$ as well as in the $L^2$-norm on $[0,1]$. In both cases, the samples are function evaluations in the Legendre nodes. The left panel is based on interpolation, the right panel corresponds to a discrete least squares approximation with $M=2N$ and regularization threshold $\epsilon = 2e^{-13}$.}\label{fig:figureB}
\end{center}
\end{figure}

In Fig.\ \ref{fig:figureB} we demonstrate the benefits of this frame for approximating the singular function given by \eqref{f_example}. Here, rather than the inner products \R{sampleONB}, we take the sampling functionals to be pointwise evaluations at the Legendre nodes, i.e., at the roots of a high degree Legendre polynomial, mapped from $[-1,1]$ to $[0,1]$.
We also compare these approximations with the polynomial interpolant at these nodes.  As is evident, polynomial interpolation performs poorly. In contrast, the convergence for the ONB$+K$ frame is significantly faster. Fig.\ \ref{fig:figureB} also illustrates the stability of the numerical approximation using oversampling with $M=2N$. The accuracy reaches machine precision, in spite of it requiring the solution of an extremely ill-conditioned linear system of equations, and this high level of accuracy is maintained as $N$ grows.

\begin{figure}[t]
\begin{center}
\subfigure[Legendre points]{
\begin{tiny}
\begin{tikzpicture}[]
\begin{axis}[height = \figureheight, ylabel = {}, xmin = {43.83}, xmax = {85.17}, ymax = {32.359365692962825}, ymode = {log}, xlabel = {M}, unbounded coords=jump,scaled x ticks = false,xlabel style = {font = {\fontsize{8 pt}{10.4 pt}\selectfont}, color = {rgb,1:red,0.00000000;green,0.00000000;blue,0.00000000}, draw opacity = 1.0, rotate = 0.0},xmajorgrids = true,xtick = {50.0,60.0,70.0,80.0},xticklabels = {$50$,$60$,$70$,$80$},xtick align = inside,xticklabel style = {font = {\fontsize{8 pt}{10.4 pt}\selectfont}, color = {rgb,1:red,0.00000000;green,0.00000000;blue,0.00000000}, draw opacity = 1.0, rotate = 0.0},x grid style = {color = {rgb,1:red,0.00000000;green,0.00000000;blue,0.00000000},
draw opacity = 0.1,
line width = 0.5,
solid},axis x line* = left,x axis line style = {color = {rgb,1:red,0.00000000;green,0.00000000;blue,0.00000000},
draw opacity = 1.0,
line width = 1,
solid},scaled y ticks = false,ylabel style = {font = {\fontsize{8 pt}{10.4 pt}\selectfont}, color = {rgb,1:red,0.00000000;green,0.00000000;blue,0.00000000}, draw opacity = 1.0, rotate = 0.0},log basis y=10,ymajorgrids = true,ytick = {1.0e-15,1.0e-10,1.0e-5,1.0},yticklabels = {$10^{-15}$,$10^{-10}$,$10^{-5}$,$10^{0}$},ytick align = inside,yticklabel style = {font = {\fontsize{8 pt}{10.4 pt}\selectfont}, color = {rgb,1:red,0.00000000;green,0.00000000;blue,0.00000000}, draw opacity = 1.0, rotate = 0.0},y grid style = {color = {rgb,1:red,0.00000000;green,0.00000000;blue,0.00000000},
draw opacity = 0.1,
line width = 0.5,
solid},axis y line* = left,y axis line style = {color = {rgb,1:red,0.00000000;green,0.00000000;blue,0.00000000},
draw opacity = 1.0,
line width = 1,
solid},    xshift = 0.0mm,
    yshift = 0.0mm,
    axis background/.style={fill={rgb,1:red,1.00000000;green,1.00000000;blue,1.00000000}}
,legend style = {color = {rgb,1:red,0.00000000;green,0.00000000;blue,0.00000000},
draw opacity = 1.0,
line width = 1,
solid,fill = {rgb,1:red,1.00000000;green,1.00000000;blue,1.00000000},fill opacity = 1.0,text opacity = 1.0,font = {\fontsize{8 pt}{10.4 pt}\selectfont}},colorbar style={title=}, ymin = {3.090295432513579e-17}, width = \figurewidth]\addplot+ [color = {rgb,1:red,0.00000000;green,0.60560316;blue,0.97868012},
draw opacity = 1.0,
line width = 2,
solid,mark = none,
mark size = 2.0,
mark options = {
            color = {rgb,1:red,0.00000000;green,0.00000000;blue,0.00000000}, draw opacity = 1.0,
            fill = {rgb,1:red,0.00000000;green,0.60560316;blue,0.97868012}, fill opacity = 1.0,
            line width = 1,
            rotate = 0,
            solid
        }]coordinates {
(45.0, 3.2307045927382205e-11)
(46.0, 8.681944052568724e-13)
(47.0, 3.19211324040225e-12)
(48.0, 8.313350008393172e-13)
(49.0, 4.676259379721159e-13)
(50.0, 2.5135449277513544e-13)
(51.0, 8.304468224196171e-14)
(52.0, 1.8740564655672642e-13)
(53.0, 9.903189379656396e-14)
(54.0, 5.240252676230739e-14)
(55.0, 1.9984014443252818e-14)
(56.0, 1.509903313490213e-14)
(57.0, 8.881784197001252e-15)
(58.0, 1.021405182655144e-14)
(59.0, 6.661338147750939e-15)
(60.0, 3.9968028886505635e-15)
(61.0, 1.3322676295501878e-15)
(62.0, 2.220446049250313e-16)
(63.0, 3.552713678800501e-15)
(64.0, 8.881784197001252e-16)
(65.0, 8.43769498715119e-15)
(66.0, 6.661338147750939e-15)
(67.0, 3.9968028886505635e-15)
(68.0, 1.3322676295501878e-14)
(69.0, 5.329070518200751e-15)
(70.0, 1.1102230246251565e-14)
(71.0, 2.353672812205332e-14)
(72.0, 2.6201263381153694e-14)
(73.0, 3.1086244689504383e-15)
(74.0, 5.773159728050814e-15)
(75.0, 2.042810365310288e-14)
(76.0, 2.4868995751603507e-14)
(77.0, 1.1546319456101628e-14)
(78.0, 4.440892098500626e-16)
(79.0, 7.105427357601002e-15)
(80.0, 8.881784197001252e-15)
(81.0, 9.769962616701378e-15)
(82.0, 8.43769498715119e-15)
(83.0, 1.6431300764452317e-14)
(84.0, 1.199040866595169e-14)
};
\addlegendentry{x=0.01}
\addlegendentry{x=0.1}
\addlegendentry{x=0.9}
\addlegendentry{L2-norm}
\addplot+ [color = {rgb,1:red,0.88887350;green,0.43564919;blue,0.27812294},
draw opacity = 1.0,
line width = 2,
solid,mark = none,
mark size = 2.0,
mark options = {
            color = {rgb,1:red,0.00000000;green,0.00000000;blue,0.00000000}, draw opacity = 1.0,
            fill = {rgb,1:red,0.88887350;green,0.43564919;blue,0.27812294}, fill opacity = 1.0,
            line width = 1,
            rotate = 0,
            solid
        }]coordinates {
(45.0, 3.6415315207705135e-14)
(46.0, 8.881784197001252e-16)
(47.0, 7.327471962526033e-15)
(48.0, 4.440892098500626e-16)
(49.0, 2.6645352591003757e-15)
(50.0, 3.3306690738754696e-15)
(51.0, 8.881784197001252e-16)
(52.0, 5.551115123125783e-15)
(53.0, 6.661338147750939e-16)
(54.0, 2.220446049250313e-15)
(55.0, 1.1102230246251565e-15)
(56.0, 1.1102230246251565e-15)
(57.0, 1.7763568394002505e-15)
(58.0, 1.7763568394002505e-15)
(59.0, 1.1102230246251565e-15)
(60.0, 1.1102230246251565e-15)
(61.0, 2.220446049250313e-16)
(62.0, 4.440892098500626e-16)
(63.0, 5.551115123125783e-15)
(64.0, 3.552713678800501e-15)
(65.0, 3.9968028886505635e-15)
(66.0, 1.7763568394002505e-15)
(67.0, 3.552713678800501e-15)
(68.0, 1.9984014443252818e-15)
(69.0, 2.6645352591003757e-15)
(70.0, 3.9968028886505635e-15)
(71.0, 3.9968028886505635e-15)
(72.0, 8.881784197001252e-15)
(73.0, 3.552713678800501e-15)
(74.0, 8.881784197001252e-16)
(75.0, 3.3306690738754696e-15)
(76.0, 6.661338147750939e-16)
(77.0, 3.552713678800501e-15)
(78.0, 5.329070518200751e-15)
(79.0, 2.220446049250313e-15)
(80.0, 1.5543122344752192e-15)
(81.0, 6.661338147750939e-16)
(82.0, 2.886579864025407e-15)
(83.0, 1.1102230246251565e-15)
(84.0, 4.440892098500626e-16)
};
\addlegendentry{x=0.01}
\addlegendentry{x=0.1}
\addlegendentry{x=0.9}
\addlegendentry{L2-norm}
\addplot+ [color = {rgb,1:red,0.24222430;green,0.64327509;blue,0.30444865},
draw opacity = 1.0,
line width = 2,
solid,mark = none,
mark size = 2.0,
mark options = {
            color = {rgb,1:red,0.00000000;green,0.00000000;blue,0.00000000}, draw opacity = 1.0,
            fill = {rgb,1:red,0.24222430;green,0.64327509;blue,0.30444865}, fill opacity = 1.0,
            line width = 1,
            rotate = 0,
            solid
        }]coordinates {
(45.0, 1.7763568394002505e-15)
(46.0, 6.217248937900877e-15)
(47.0, 8.881784197001252e-16)
(48.0, 5.773159728050814e-15)
(49.0, 6.217248937900877e-15)
(50.0, 1.7763568394002505e-15)
(51.0, 8.881784197001252e-16)
(52.0, 2.220446049250313e-16)
(53.0, 2.6645352591003757e-15)
(54.0, 8.881784197001252e-16)
(55.0, 3.9968028886505635e-15)
(56.0, 2.6645352591003757e-15)
(57.0, 3.9968028886505635e-15)
(58.0, 4.440892098500626e-15)
(59.0, 1.7763568394002505e-15)
(60.0, 1.7763568394002505e-15)
(61.0, 2.6645352591003757e-15)
(62.0, 3.552713678800501e-15)
(63.0, 1.7763568394002505e-15)
(64.0, 4.440892098500626e-16)
(65.0, 4.440892098500626e-16)
(66.0, 3.1086244689504383e-15)
(67.0, 2.6645352591003757e-15)
(68.0, 4.440892098500626e-15)
(69.0, 8.881784197001252e-16)
(70.0, 1.7763568394002505e-15)
(71.0, 1.3322676295501878e-15)
(72.0, 8.881784197001252e-16)
(73.0, 8.881784197001252e-16)
(74.0, 4.440892098500626e-16)
(75.0, 1.3322676295501878e-15)
(76.0, 2.220446049250313e-15)
(77.0, 2.220446049250313e-16)
(78.0, 1.3322676295501878e-15)
(79.0, 2.6645352591003757e-15)
(80.0, 1.7763568394002505e-15)
(81.0, 4.440892098500626e-16)
(82.0, 2.6645352591003757e-15)
(83.0, 4.440892098500626e-16)
(84.0, 4.440892098500626e-15)
};
\addlegendentry{x=0.01}
\addlegendentry{x=0.1}
\addlegendentry{x=0.9}
\addlegendentry{L2-norm}
\addplot+ [color = {rgb,1:red,0.76444018;green,0.44411178;blue,0.82429754},
draw opacity = 1.0,
line width = 2,
solid,mark = none,
mark size = 2.0,
mark options = {
            color = {rgb,1:red,0.00000000;green,0.00000000;blue,0.00000000}, draw opacity = 1.0,
            fill = {rgb,1:red,0.76444018;green,0.44411178;blue,0.82429754}, fill opacity = 1.0,
            line width = 1,
            rotate = 0,
            solid
        }]coordinates {
(45.0, 4.2820295569880965e-9)
(46.0, 7.48160734281687e-11)
(47.0, 2.877877033388777e-10)
(48.0, 7.148039999453229e-11)
(49.0, 4.3698089759436547e-11)
(50.0, 2.2655779441610173e-11)
(51.0, 4.3028361743343875e-12)
(52.0, 1.3902207180704202e-11)
(53.0, 1.0619290970864144e-11)
(54.0, 4.767036593612669e-12)
(55.0, 3.5911264995196155e-12)
(56.0, 3.94877191840535e-12)
(57.0, 3.9314344972375396e-12)
(58.0, 3.893919425511561e-12)
(59.0, 2.5509279794639403e-12)
(60.0, 2.3664352573073574e-12)
(61.0, 1.267839593292652e-12)
(62.0, 2.0716721680023443e-12)
(63.0, 1.6259685945095406e-12)
(64.0, 1.051557212208938e-12)
(65.0, 1.2627180522388045e-12)
(66.0, 1.0125923106801077e-12)
(67.0, 1.1671209157167524e-12)
(68.0, 1.1407276176861498e-12)
(69.0, 7.728646536219936e-13)
(70.0, 6.509976674094684e-13)
(71.0, 7.831754676431424e-13)
(72.0, 5.736876777668535e-13)
(73.0, 4.93010514969452e-13)
(74.0, 4.978470380053252e-13)
(75.0, 3.645041019742801e-13)
(76.0, 2.864297393324223e-13)
(77.0, 1.5838383571242735e-13)
(78.0, 2.7319712560377096e-13)
(79.0, 1.787814876923848e-13)
(80.0, 1.814241987866066e-13)
(81.0, 1.6636982546683266e-13)
(82.0, 1.7689932278442484e-14)
(83.0, 5.944422863380827e-14)
(84.0, 7.228540329396309e-15)
};
\addlegendentry{x=0.01}
\addlegendentry{x=0.1}
\addlegendentry{x=0.9}
\addlegendentry{L2-norm}
\end{axis}

\end{tikzpicture}
\end{tiny}
}
\subfigure[Equispaced points]{
\begin{tiny}
\input{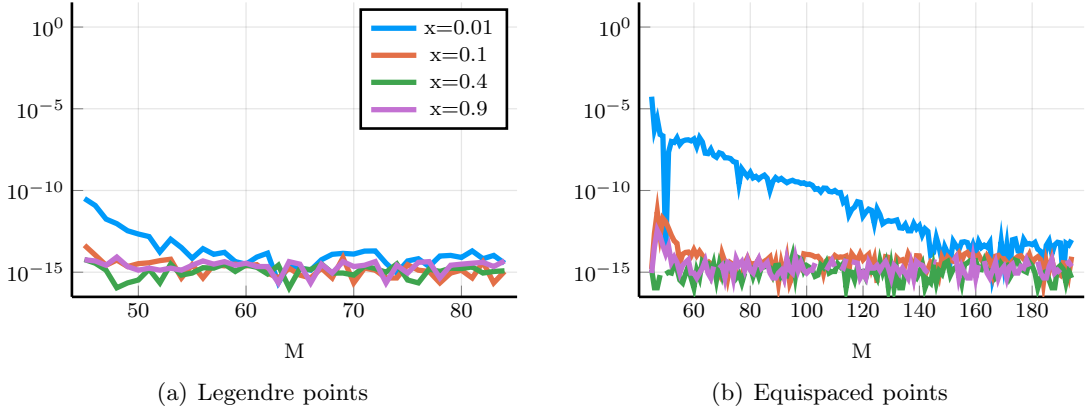}
\end{tiny}
}
\caption{Pointwise error as a function of the number of samples $M$ for the ONB$+K$ frame based on Legendre polynomials, with degree $N=40$, $K=5$ and $w(x)=\log x$ as in Fig.\ \ref{fig:figureB}(b). Similar to regular polynomial approximation, using Legendre points yields better accuracy than using equispaced points. They also require less oversampling, i.e., smaller values of $M$, to achieve the best error.}\label{fig:oversampling}
\end{center}
\end{figure}

The influence of the oversampling factor $M$ is illustrated in Fig.\ \ref{fig:oversampling}. Here, the error is shown as a function of $M$, for constant $N=40$. Best accuracy is only achieved for $M > N$, i.e., when using some amount of oversampling. We have used discrete sampling in this figure using Legendre points (left panel) and equispaced points on $[0,1]$ (right panel). It is not unexpected that Legendre points are a better choice: less oversampling is needed to achieve the best accuracy for the given $N$.

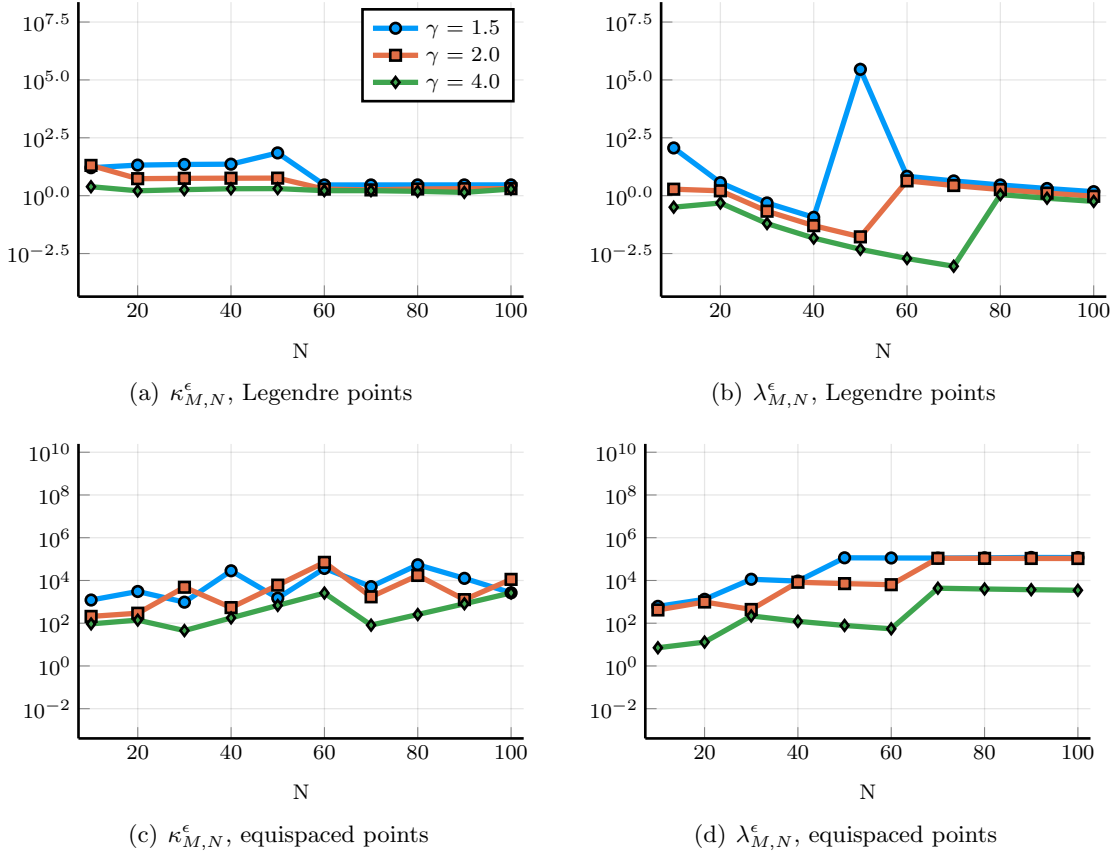
\begin{figure}[t]
\begin{center}
\subfigure[$\protect \kappa_{M,N}^\epsilon$, Legendre points]{
\begin{tiny}
\begin{tikzpicture}[]
\begin{axis}[height = \figureheight, ylabel = {}, xmin = {7.300000000000001}, xmax = {102.7}, ymax = {2.29086765276777e8}, ymode = {log}, xlabel = {N}, unbounded coords=jump,scaled x ticks = false,xlabel style = {font = {\fontsize{8 pt}{10.4 pt}\selectfont}, color = {rgb,1:red,0.00000000;green,0.00000000;blue,0.00000000}, draw opacity = 1.0, rotate = 0.0},xmajorgrids = true,xtick = {20.0,40.0,60.0,80.0,100.0},xticklabels = {$20$,$40$,$60$,$80$,$100$},xtick align = inside,xticklabel style = {font = {\fontsize{8 pt}{10.4 pt}\selectfont}, color = {rgb,1:red,0.00000000;green,0.00000000;blue,0.00000000}, draw opacity = 1.0, rotate = 0.0},x grid style = {color = {rgb,1:red,0.00000000;green,0.00000000;blue,0.00000000},
draw opacity = 0.1,
line width = 0.5,
solid},axis lines* = left,x axis line style = {color = {rgb,1:red,0.00000000;green,0.00000000;blue,0.00000000},
draw opacity = 1.0,
line width = 1,
solid},scaled y ticks = false,ylabel style = {font = {\fontsize{8 pt}{10.4 pt}\selectfont}, color = {rgb,1:red,0.00000000;green,0.00000000;blue,0.00000000}, draw opacity = 1.0, rotate = 0.0},log basis y=10,ymajorgrids = true,ytick = {0.0031622776601683794,1.0,316.22776601683796,100000.0,3.162277660168379e7},yticklabels = {$10^{-2.5}$,$10^{0.0}$,$10^{2.5}$,$10^{5.0}$,$10^{7.5}$},ytick align = inside,yticklabel style = {font = {\fontsize{8 pt}{10.4 pt}\selectfont}, color = {rgb,1:red,0.00000000;green,0.00000000;blue,0.00000000}, draw opacity = 1.0, rotate = 0.0},y grid style = {color = {rgb,1:red,0.00000000;green,0.00000000;blue,0.00000000},
draw opacity = 0.1,
line width = 0.5,
solid},axis lines* = left,y axis line style = {color = {rgb,1:red,0.00000000;green,0.00000000;blue,0.00000000},
draw opacity = 1.0,
line width = 1,
solid},    xshift = 0.0mm,
    yshift = 0.0mm,
    axis background/.style={fill={rgb,1:red,1.00000000;green,1.00000000;blue,1.00000000}}
,legend style = {color = {rgb,1:red,0.00000000;green,0.00000000;blue,0.00000000},
draw opacity = 1.0,
line width = 1,
solid,fill = {rgb,1:red,1.00000000;green,1.00000000;blue,1.00000000},font = {\fontsize{8 pt}{10.4 pt}\selectfont}},colorbar style={title=}, ymin = {4.3651583224016566e-5}, width = \figurewidth]\addplot+ [color = {rgb,1:red,0.00000000;green,0.60560316;blue,0.97868012},
draw opacity = 1.0,
line width = 2,
solid,mark = *,
mark size = 2.0,
mark options = {
    color = {rgb,1:red,0.00000000;green,0.00000000;blue,0.00000000}, draw opacity = 1.0,
    fill = {rgb,1:red,0.00000000;green,0.60560316;blue,0.97868012}, fill opacity = 1.0,
    line width = 1,
    rotate = 0,
    solid
}]coordinates {
(10.0, 16.155606618294794)
(20.0, 20.95182141027372)
(30.0, 22.31857932013877)
(40.0, 22.949893782847653)
(50.0, 70.43322040156151)
(60.0, 2.9006961454653934)
(70.0, 2.911982822134397)
(80.0, 2.9173842267285597)
(90.0, 2.918775571457945)
(100.0, 2.9242750313113066)
};
\addlegendentry{$\gamma$ = 1.5}
\addplot+ [color = {rgb,1:red,0.88887350;green,0.43564919;blue,0.27812294},
draw opacity = 1.0,
line width = 2,
solid,mark = square*,
mark size = 2.0,
mark options = {
    color = {rgb,1:red,0.00000000;green,0.00000000;blue,0.00000000}, draw opacity = 1.0,
    fill = {rgb,1:red,0.88887350;green,0.43564919;blue,0.27812294}, fill opacity = 1.0,
    line width = 1,
    rotate = 0,
    solid
}]coordinates {
(10.0, 20.295044069958706)
(20.0, 5.43362452325313)
(30.0, 5.613620766174688)
(40.0, 5.694146408024371)
(50.0, 5.735632565260178)
(60.0, 1.9210802345483098)
(70.0, 1.780125552790867)
(80.0, 1.9508825360319384)
(90.0, 2.0096538484355557)
(100.0, 2.01748228085298)
};
\addlegendentry{$\gamma$ = 2.0}
\addplot+ [color = {rgb,1:red,0.24222430;green,0.64327509;blue,0.30444865},
draw opacity = 1.0,
line width = 2,
solid,mark = diamond*,
mark size = 2.0,
mark options = {
    color = {rgb,1:red,0.00000000;green,0.00000000;blue,0.00000000}, draw opacity = 1.0,
    fill = {rgb,1:red,0.24222430;green,0.64327509;blue,0.30444865}, fill opacity = 1.0,
    line width = 1,
    rotate = 0,
    solid
}]coordinates {
(10.0, 2.4450429641001805)
(20.0, 1.6291115728328072)
(30.0, 1.8247758984603055)
(40.0, 1.9971297390247733)
(50.0, 2.022006694897507)
(60.0, 1.6486749770578963)
(70.0, 1.6502349859752206)
(80.0, 1.545999119811728)
(90.0, 1.3897964591498517)
(100.0, 1.934512626468712)
};
\addlegendentry{$\gamma$ = 4.0}
\end{axis}

\end{tikzpicture}
\end{tiny}
}
\subfigure[$\protect \lambda_{M,N}^\epsilon$, Legendre points]{
\begin{tiny}
\begin{tikzpicture}[]
\begin{axis}[height = \figureheight, ylabel = {}, xmin = {7.300000000000001}, xmax = {102.7}, ymax = {2.29086765276777e8}, ymode = {log}, xlabel = {N}, unbounded coords=jump,scaled x ticks = false,xlabel style = {font = {\fontsize{8 pt}{10.4 pt}\selectfont}, color = {rgb,1:red,0.00000000;green,0.00000000;blue,0.00000000}, draw opacity = 1.0, rotate = 0.0},xmajorgrids = true,xtick = {20.0,40.0,60.0,80.0,100.0},xticklabels = {$20$,$40$,$60$,$80$,$100$},xtick align = inside,xticklabel style = {font = {\fontsize{8 pt}{10.4 pt}\selectfont}, color = {rgb,1:red,0.00000000;green,0.00000000;blue,0.00000000}, draw opacity = 1.0, rotate = 0.0},x grid style = {color = {rgb,1:red,0.00000000;green,0.00000000;blue,0.00000000},
draw opacity = 0.1,
line width = 0.5,
solid},axis lines* = left,x axis line style = {color = {rgb,1:red,0.00000000;green,0.00000000;blue,0.00000000},
draw opacity = 1.0,
line width = 1,
solid},scaled y ticks = false,ylabel style = {font = {\fontsize{8 pt}{10.4 pt}\selectfont}, color = {rgb,1:red,0.00000000;green,0.00000000;blue,0.00000000}, draw opacity = 1.0, rotate = 0.0},log basis y=10,ymajorgrids = true,ytick = {0.0031622776601683794,1.0,316.22776601683796,100000.0,3.162277660168379e7},yticklabels = {$10^{-2.5}$,$10^{0.0}$,$10^{2.5}$,$10^{5.0}$,$10^{7.5}$},ytick align = inside,yticklabel style = {font = {\fontsize{8 pt}{10.4 pt}\selectfont}, color = {rgb,1:red,0.00000000;green,0.00000000;blue,0.00000000}, draw opacity = 1.0, rotate = 0.0},y grid style = {color = {rgb,1:red,0.00000000;green,0.00000000;blue,0.00000000},
draw opacity = 0.1,
line width = 0.5,
solid},axis lines* = left,y axis line style = {color = {rgb,1:red,0.00000000;green,0.00000000;blue,0.00000000},
draw opacity = 1.0,
line width = 1,
solid},    xshift = 0.0mm,
    yshift = 0.0mm,
    axis background/.style={fill={rgb,1:red,1.00000000;green,1.00000000;blue,1.00000000}}
,legend style = {color = {rgb,1:red,0.00000000;green,0.00000000;blue,0.00000000},
draw opacity = 1.0,
line width = 1,
solid,fill = {rgb,1:red,1.00000000;green,1.00000000;blue,1.00000000},font = {\fontsize{8 pt}{10.4 pt}\selectfont}},colorbar style={title=}, ymin = {4.3651583224016566e-5}, width = \figurewidth]\addplot+ [color = {rgb,1:red,0.00000000;green,0.60560316;blue,0.97868012},
draw opacity = 1.0,
line width = 2,
solid,mark = *,
mark size = 2.0,
mark options = {
    color = {rgb,1:red,0.00000000;green,0.00000000;blue,0.00000000}, draw opacity = 1.0,
    fill = {rgb,1:red,0.00000000;green,0.60560316;blue,0.97868012}, fill opacity = 1.0,
    line width = 1,
    rotate = 0,
    solid
},forget plot]coordinates {
(10.0, 115.3079786736401)
(20.0, 3.6576161912277843)
(30.0, 0.4864944990536276)
(40.0, 0.11616275846012589)
(50.0, 286204.7123944765)
(60.0, 6.877936947944428)
(70.0, 4.335289689094637)
(80.0, 2.906480311089972)
(90.0, 2.0425723254322503)
(100.0, 1.4898014288526324)
};
\addplot+ [color = {rgb,1:red,0.88887350;green,0.43564919;blue,0.27812294},
draw opacity = 1.0,
line width = 2,
solid,mark = square*,
mark size = 2.0,
mark options = {
    color = {rgb,1:red,0.00000000;green,0.00000000;blue,0.00000000}, draw opacity = 1.0,
    fill = {rgb,1:red,0.88887350;green,0.43564919;blue,0.27812294}, fill opacity = 1.0,
    line width = 1,
    rotate = 0,
    solid
},forget plot]coordinates {
(10.0, 1.9071728009659081)
(20.0, 1.6243820753665505)
(30.0, 0.2138180664009403)
(40.0, 0.05081404109710444)
(50.0, 0.016674261977347705)
(60.0, 4.3671899409476165)
(70.0, 2.7504658562144093)
(80.0, 1.8428670899208581)
(90.0, 1.2945050997769831)
(100.0, 0.9438356538907806)
};
\addplot+ [color = {rgb,1:red,0.24222430;green,0.64327509;blue,0.30444865},
draw opacity = 1.0,
line width = 2,
solid,mark = diamond*,
mark size = 2.0,
mark options = {
    color = {rgb,1:red,0.00000000;green,0.00000000;blue,0.00000000}, draw opacity = 1.0,
    fill = {rgb,1:red,0.24222430;green,0.64327509;blue,0.30444865}, fill opacity = 1.0,
    line width = 1,
    rotate = 0,
    solid
},forget plot]coordinates {
(10.0, 0.3186420534121448)
(20.0, 0.48435073819613667)
(30.0, 0.0626245284492879)
(40.0, 0.014771766262842411)
(50.0, 0.004827691026625592)
(60.0, 0.001937529613820747)
(70.0, 0.0008957421011521329)
(80.0, 1.1083121609745987)
(90.0, 0.7780700070971203)
(100.0, 0.5670419726539707)
};
\end{axis}

\end{tikzpicture}
\end{tiny}
}
\vspace{0.5cm}
\subfigure[$\protect \kappa_{M,N}^\epsilon$, equispaced points]{
\begin{tiny}
\begin{tikzpicture}[]
\begin{axis}[height = \figureheight, ylabel = {}, xmin = {7.300000000000001}, xmax = {102.7}, ymax = {2.4547089156850338e10}, ymode = {log}, xlabel = {N}, unbounded coords=jump,scaled x ticks = false,xlabel style = {font = {\fontsize{8 pt}{10.4 pt}\selectfont}, color = {rgb,1:red,0.00000000;green,0.00000000;blue,0.00000000}, draw opacity = 1.0, rotate = 0.0},xmajorgrids = true,xtick = {20.0,40.0,60.0,80.0,100.0},xticklabels = {$20$,$40$,$60$,$80$,$100$},xtick align = inside,xticklabel style = {font = {\fontsize{8 pt}{10.4 pt}\selectfont}, color = {rgb,1:red,0.00000000;green,0.00000000;blue,0.00000000}, draw opacity = 1.0, rotate = 0.0},x grid style = {color = {rgb,1:red,0.00000000;green,0.00000000;blue,0.00000000},
draw opacity = 0.1,
line width = 0.5,
solid},axis lines* = left,x axis line style = {color = {rgb,1:red,0.00000000;green,0.00000000;blue,0.00000000},
draw opacity = 1.0,
line width = 1,
solid},scaled y ticks = false,ylabel style = {font = {\fontsize{8 pt}{10.4 pt}\selectfont}, color = {rgb,1:red,0.00000000;green,0.00000000;blue,0.00000000}, draw opacity = 1.0, rotate = 0.0},log basis y=10,ymajorgrids = true,ytick = {0.01,1.0,100.0,10000.0,1.0e6,1.0e8,1.0e10},yticklabels = {$10^{-2}$,$10^{0}$,$10^{2}$,$10^{4}$,$10^{6}$,$10^{8}$,$10^{10}$},ytick align = inside,yticklabel style = {font = {\fontsize{8 pt}{10.4 pt}\selectfont}, color = {rgb,1:red,0.00000000;green,0.00000000;blue,0.00000000}, draw opacity = 1.0, rotate = 0.0},y grid style = {color = {rgb,1:red,0.00000000;green,0.00000000;blue,0.00000000},
draw opacity = 0.1,
line width = 0.5,
solid},axis lines* = left,y axis line style = {color = {rgb,1:red,0.00000000;green,0.00000000;blue,0.00000000},
draw opacity = 1.0,
line width = 1,
solid},    xshift = 0.0mm,
    yshift = 0.0mm,
    axis background/.style={fill={rgb,1:red,1.00000000;green,1.00000000;blue,1.00000000}}
,legend style = {color = {rgb,1:red,0.00000000;green,0.00000000;blue,0.00000000},
draw opacity = 1.0,
line width = 1,
solid,fill = {rgb,1:red,1.00000000;green,1.00000000;blue,1.00000000},font = {\fontsize{8 pt}{10.4 pt}\selectfont}},colorbar style={title=}, ymin = {0.0004073802778041126}, width = \figurewidth]\addplot+ [color = {rgb,1:red,0.00000000;green,0.60560316;blue,0.97868012},
draw opacity = 1.0,
line width = 2,
solid,mark = *,
mark size = 2.0,
mark options = {
    color = {rgb,1:red,0.00000000;green,0.00000000;blue,0.00000000}, draw opacity = 1.0,
    fill = {rgb,1:red,0.00000000;green,0.60560316;blue,0.97868012}, fill opacity = 1.0,
    line width = 1,
    rotate = 0,
    solid
},forget plot]coordinates {
(10.0, 1216.3366362167978)
(20.0, 3053.1524777902423)
(30.0, 965.5785617812206)
(40.0, 28378.703533918037)
(50.0, 1492.3371532068388)
(60.0, 36690.01363872824)
(70.0, 5147.829826629525)
(80.0, 54314.150984603344)
(90.0, 12560.788155470407)
(100.0, 2651.7631834794547)
};
\addplot+ [color = {rgb,1:red,0.88887350;green,0.43564919;blue,0.27812294},
draw opacity = 1.0,
line width = 2,
solid,mark = square*,
mark size = 2.0,
mark options = {
    color = {rgb,1:red,0.00000000;green,0.00000000;blue,0.00000000}, draw opacity = 1.0,
    fill = {rgb,1:red,0.88887350;green,0.43564919;blue,0.27812294}, fill opacity = 1.0,
    line width = 1,
    rotate = 0,
    solid
},forget plot]coordinates {
(10.0, 208.39086033661562)
(20.0, 293.9091114299237)
(30.0, 4824.449661583019)
(40.0, 538.7104827213077)
(50.0, 6043.129836019676)
(60.0, 70498.95231341937)
(70.0, 1717.9167533134355)
(80.0, 17193.249442426455)
(90.0, 1269.742374523129)
(100.0, 11254.886385407082)
};
\addplot+ [color = {rgb,1:red,0.24222430;green,0.64327509;blue,0.30444865},
draw opacity = 1.0,
line width = 2,
solid,mark = diamond*,
mark size = 2.0,
mark options = {
    color = {rgb,1:red,0.00000000;green,0.00000000;blue,0.00000000}, draw opacity = 1.0,
    fill = {rgb,1:red,0.24222430;green,0.64327509;blue,0.30444865}, fill opacity = 1.0,
    line width = 1,
    rotate = 0,
    solid
},forget plot]coordinates {
(10.0, 92.95855855018357)
(20.0, 140.13762887784222)
(30.0, 45.12794809722437)
(40.0, 177.46414259769674)
(50.0, 680.8683750327995)
(60.0, 2572.6317620593404)
(70.0, 81.30511869862237)
(80.0, 254.98727833837648)
(90.0, 812.3189270332504)
(100.0, 2619.082718712736)
};
\end{axis}

\end{tikzpicture}
\end{tiny}
}
\subfigure[$\protect \lambda_{M,N}^\epsilon$, equispaced points]{
\begin{tiny}
\begin{tikzpicture}[]
\begin{axis}[height = \figureheight, ylabel = {}, xmin = {7.300000000000001}, xmax = {102.7}, ymax = {2.4547089156850338e10}, ymode = {log}, xlabel = {N}, unbounded coords=jump,scaled x ticks = false,xlabel style = {font = {\fontsize{8 pt}{10.4 pt}\selectfont}, color = {rgb,1:red,0.00000000;green,0.00000000;blue,0.00000000}, draw opacity = 1.0, rotate = 0.0},xmajorgrids = true,xtick = {20.0,40.0,60.0,80.0,100.0},xticklabels = {$20$,$40$,$60$,$80$,$100$},xtick align = inside,xticklabel style = {font = {\fontsize{8 pt}{10.4 pt}\selectfont}, color = {rgb,1:red,0.00000000;green,0.00000000;blue,0.00000000}, draw opacity = 1.0, rotate = 0.0},x grid style = {color = {rgb,1:red,0.00000000;green,0.00000000;blue,0.00000000},
draw opacity = 0.1,
line width = 0.5,
solid},axis lines* = left,x axis line style = {color = {rgb,1:red,0.00000000;green,0.00000000;blue,0.00000000},
draw opacity = 1.0,
line width = 1,
solid},scaled y ticks = false,ylabel style = {font = {\fontsize{8 pt}{10.4 pt}\selectfont}, color = {rgb,1:red,0.00000000;green,0.00000000;blue,0.00000000}, draw opacity = 1.0, rotate = 0.0},log basis y=10,ymajorgrids = true,ytick = {0.01,1.0,100.0,10000.0,1.0e6,1.0e8,1.0e10},yticklabels = {$10^{-2}$,$10^{0}$,$10^{2}$,$10^{4}$,$10^{6}$,$10^{8}$,$10^{10}$},ytick align = inside,yticklabel style = {font = {\fontsize{8 pt}{10.4 pt}\selectfont}, color = {rgb,1:red,0.00000000;green,0.00000000;blue,0.00000000}, draw opacity = 1.0, rotate = 0.0},y grid style = {color = {rgb,1:red,0.00000000;green,0.00000000;blue,0.00000000},
draw opacity = 0.1,
line width = 0.5,
solid},axis lines* = left,y axis line style = {color = {rgb,1:red,0.00000000;green,0.00000000;blue,0.00000000},
draw opacity = 1.0,
line width = 1,
solid},    xshift = 0.0mm,
    yshift = 0.0mm,
    axis background/.style={fill={rgb,1:red,1.00000000;green,1.00000000;blue,1.00000000}}
,legend style = {color = {rgb,1:red,0.00000000;green,0.00000000;blue,0.00000000},
draw opacity = 1.0,
line width = 1,
solid,fill = {rgb,1:red,1.00000000;green,1.00000000;blue,1.00000000},font = {\fontsize{8 pt}{10.4 pt}\selectfont}},colorbar style={title=}, ymin = {0.0004073802778041126}, width = \figurewidth]\addplot+ [color = {rgb,1:red,0.00000000;green,0.60560316;blue,0.97868012},
draw opacity = 1.0,
line width = 2,
solid,mark = *,
mark size = 2.0,
mark options = {
    color = {rgb,1:red,0.00000000;green,0.00000000;blue,0.00000000}, draw opacity = 1.0,
    fill = {rgb,1:red,0.00000000;green,0.60560316;blue,0.97868012}, fill opacity = 1.0,
    line width = 1,
    rotate = 0,
    solid
},forget plot]coordinates {
(10.0, 614.4182890416579)
(20.0, 1312.0204183764292)
(30.0, 11386.806002159497)
(40.0, 9344.17149619305)
(50.0, 115044.03614834679)
(60.0, 112938.71834161373)
(70.0, 111386.77579744397)
(80.0, 114866.28695648452)
(90.0, 119280.90100913192)
(100.0, 117937.14226292179)
};
\addplot+ [color = {rgb,1:red,0.88887350;green,0.43564919;blue,0.27812294},
draw opacity = 1.0,
line width = 2,
solid,mark = square*,
mark size = 2.0,
mark options = {
    color = {rgb,1:red,0.00000000;green,0.00000000;blue,0.00000000}, draw opacity = 1.0,
    fill = {rgb,1:red,0.88887350;green,0.43564919;blue,0.27812294}, fill opacity = 1.0,
    line width = 1,
    rotate = 0,
    solid
},forget plot]coordinates {
(10.0, 417.14479893916695)
(20.0, 983.6722070632029)
(30.0, 439.2433923118027)
(40.0, 8255.395138597522)
(50.0, 7131.230948747635)
(60.0, 6344.726477179758)
(70.0, 108948.2983224817)
(80.0, 108038.68895734684)
(90.0, 107235.93262841708)
(100.0, 106708.15956070744)
};
\addplot+ [color = {rgb,1:red,0.24222430;green,0.64327509;blue,0.30444865},
draw opacity = 1.0,
line width = 2,
solid,mark = diamond*,
mark size = 2.0,
mark options = {
    color = {rgb,1:red,0.00000000;green,0.00000000;blue,0.00000000}, draw opacity = 1.0,
    fill = {rgb,1:red,0.24222430;green,0.64327509;blue,0.30444865}, fill opacity = 1.0,
    line width = 1,
    rotate = 0,
    solid
},forget plot]coordinates {
(10.0, 7.104850899906832)
(20.0, 13.094218246125257)
(30.0, 219.60953273646322)
(40.0, 122.47512077667783)
(50.0, 78.48957149969718)
(60.0, 54.82566225448668)
(70.0, 4339.1466162111155)
(80.0, 3989.3908979319117)
(90.0, 3706.7808993299745)
(100.0, 3472.6261772295784)
};
\end{axis}

\end{tikzpicture}
\end{tiny}
}
\caption{The values of $\kappa_{M,N}^\epsilon$ and $\lambda_{M,N}^\epsilon$ are shown as a function of $N$ with constant oversampling $M = \gamma N$ and varying factors $\gamma$. Legendre points (top row) and equispaced points (bottom row) are used for the ONB$+K$ frame using Legendre polynomials, $w(x)=\log(x)$ and $K=5$. The threshold used here is $\epsilon = 1e^{-5}$: the values are bounded for Legendre points but they approach $1/\epsilon$ for equispaced points. Equispaced points require more than linear oversampling.}\label{fig:samplingrate}
\end{center}
\end{figure}

The behaviour shown in Fig.\ \ref{fig:oversampling} can be explained by computing the corresponding constants $\kappa_{M,N}^\epsilon$ and $\lambda_{M,N}^\epsilon$. This also serves to illustrate the stable sampling rate for this problem. Their values are shown in Fig.\ \ref{fig:samplingrate} for several choices of the oversampling factor $\gamma$, with $M$ and $N$ such that $M = \gamma N$. The convergence of both values to constants of modest size, in particular much smaller than $1/\epsilon$, suggests that the stable sampling rate is indeed linear when sampling in the Legendre points. For equispaced points, as could be expected, this does not seem to be the case.

\begin{figure}[t]
\begin{center}
\subfigure[$\protect \kappa_{M,N}^\epsilon$, Legendre points]{
\begin{tiny}
\begin{tikzpicture}[]
\begin{axis}[height = \figureheight, ylabel = {}, xmin = {7.300000000000001}, xmax = {102.7}, ymax = {2.6302679918953815e10}, ymode = {log}, xlabel = {N}, unbounded coords=jump,scaled x ticks = false,xlabel style = {font = {\fontsize{8 pt}{10.4 pt}\selectfont}, color = {rgb,1:red,0.00000000;green,0.00000000;blue,0.00000000}, draw opacity = 1.0, rotate = 0.0},xmajorgrids = true,xtick = {20.0,40.0,60.0,80.0,100.0},xticklabels = {$20$,$40$,$60$,$80$,$100$},xtick align = inside,xticklabel style = {font = {\fontsize{8 pt}{10.4 pt}\selectfont}, color = {rgb,1:red,0.00000000;green,0.00000000;blue,0.00000000}, draw opacity = 1.0, rotate = 0.0},x grid style = {color = {rgb,1:red,0.00000000;green,0.00000000;blue,0.00000000},
draw opacity = 0.1,
line width = 0.5,
solid},axis lines* = left,x axis line style = {color = {rgb,1:red,0.00000000;green,0.00000000;blue,0.00000000},
draw opacity = 1.0,
line width = 1,
solid},scaled y ticks = false,ylabel style = {font = {\fontsize{8 pt}{10.4 pt}\selectfont}, color = {rgb,1:red,0.00000000;green,0.00000000;blue,0.00000000}, draw opacity = 1.0, rotate = 0.0},log basis y=10,ymajorgrids = true,ytick = {0.0031622776601683794,1.0,316.22776601683796,100000.0,3.162277660168379e7,1.0e10},yticklabels = {$10^{-2.5}$,$10^{0.0}$,$10^{2.5}$,$10^{5.0}$,$10^{7.5}$,$10^{10.0}$},ytick align = inside,yticklabel style = {font = {\fontsize{8 pt}{10.4 pt}\selectfont}, color = {rgb,1:red,0.00000000;green,0.00000000;blue,0.00000000}, draw opacity = 1.0, rotate = 0.0},y grid style = {color = {rgb,1:red,0.00000000;green,0.00000000;blue,0.00000000},
draw opacity = 0.1,
line width = 0.5,
solid},axis lines* = left,y axis line style = {color = {rgb,1:red,0.00000000;green,0.00000000;blue,0.00000000},
draw opacity = 1.0,
line width = 1,
solid},    xshift = 0.0mm,
    yshift = 0.0mm,
    axis background/.style={fill={rgb,1:red,1.00000000;green,1.00000000;blue,1.00000000}}
,legend style = {color = {rgb,1:red,0.00000000;green,0.00000000;blue,0.00000000},
draw opacity = 1.0,
line width = 1,
solid,fill = {rgb,1:red,1.00000000;green,1.00000000;blue,1.00000000},font = {\fontsize{8 pt}{10.4 pt}\selectfont}},colorbar style={title=}, ymin = {3.8018939632056124e-5}, width = \figurewidth]\addplot+ [color = {rgb,1:red,0.00000000;green,0.60560316;blue,0.97868012},
draw opacity = 1.0,
line width = 2,
solid,mark = *,
mark size = 2.0,
mark options = {
    color = {rgb,1:red,0.00000000;green,0.00000000;blue,0.00000000}, draw opacity = 1.0,
    fill = {rgb,1:red,0.00000000;green,0.60560316;blue,0.97868012}, fill opacity = 1.0,
    line width = 1,
    rotate = 0,
    solid
}]coordinates {
(10.0, 3275.2651647129314)
(20.0, 230.0264132148907)
(30.0, 245.65015778921583)
(40.0, 22.949893782847653)
(50.0, 272174.0367779492)
(60.0, 23.542303034311917)
(70.0, 23.70467929881577)
(80.0, 23.82399946429935)
(90.0, 23.91523889760706)
(100.0, 23.98780743983285)
};
\addlegendentry{$\gamma$ = 1.5}
\addplot+ [color = {rgb,1:red,0.88887350;green,0.43564919;blue,0.27812294},
draw opacity = 1.0,
line width = 2,
solid,mark = square*,
mark size = 2.0,
mark options = {
    color = {rgb,1:red,0.00000000;green,0.00000000;blue,0.00000000}, draw opacity = 1.0,
    fill = {rgb,1:red,0.88887350;green,0.43564919;blue,0.27812294}, fill opacity = 1.0,
    line width = 1,
    rotate = 0,
    solid
}]coordinates {
(10.0, 111.75160154124427)
(20.0, 24.074827469872442)
(30.0, 25.106815876740892)
(40.0, 25.57187190664811)
(50.0, 5.735632565260178)
(60.0, 5.760958266226215)
(70.0, 5.778233518403134)
(80.0, 5.791232016855182)
(90.0, 5.80068752179624)
(100.0, 5.807973654185573)
};
\addlegendentry{$\gamma$ = 2.0}
\addplot+ [color = {rgb,1:red,0.24222430;green,0.64327509;blue,0.30444865},
draw opacity = 1.0,
line width = 2,
solid,mark = diamond*,
mark size = 2.0,
mark options = {
    color = {rgb,1:red,0.00000000;green,0.00000000;blue,0.00000000}, draw opacity = 1.0,
    fill = {rgb,1:red,0.24222430;green,0.64327509;blue,0.30444865}, fill opacity = 1.0,
    line width = 1,
    rotate = 0,
    solid
}]coordinates {
(10.0, 4.405612526967264)
(20.0, 2.585634147270105)
(30.0, 2.619336146590945)
(40.0, 2.6314399852815744)
(50.0, 2.636990284274825)
(60.0, 1.6486749770578963)
(70.0, 1.6502349859752206)
(80.0, 1.6512762684704394)
(90.0, 1.651239926202366)
(100.0, 1.9350674594285246)
};
\addlegendentry{$\gamma$ = 4.0}
\end{axis}

\end{tikzpicture}
\end{tiny}
}
\subfigure[$\protect \lambda_{M,N}^\epsilon$, Legendre points]{
\begin{tiny}
\begin{tikzpicture}[]
\begin{axis}[height = \figureheight, ylabel = {}, xmin = {7.300000000000001}, xmax = {102.7}, ymax = {2.6302679918953815e10}, ymode = {log}, xlabel = {N}, unbounded coords=jump,scaled x ticks = false,xlabel style = {font = {\fontsize{8 pt}{10.4 pt}\selectfont}, color = {rgb,1:red,0.00000000;green,0.00000000;blue,0.00000000}, draw opacity = 1.0, rotate = 0.0},xmajorgrids = true,xtick = {20.0,40.0,60.0,80.0,100.0},xticklabels = {$20$,$40$,$60$,$80$,$100$},xtick align = inside,xticklabel style = {font = {\fontsize{8 pt}{10.4 pt}\selectfont}, color = {rgb,1:red,0.00000000;green,0.00000000;blue,0.00000000}, draw opacity = 1.0, rotate = 0.0},x grid style = {color = {rgb,1:red,0.00000000;green,0.00000000;blue,0.00000000},
draw opacity = 0.1,
line width = 0.5,
solid},axis lines* = left,x axis line style = {color = {rgb,1:red,0.00000000;green,0.00000000;blue,0.00000000},
draw opacity = 1.0,
line width = 1,
solid},scaled y ticks = false,ylabel style = {font = {\fontsize{8 pt}{10.4 pt}\selectfont}, color = {rgb,1:red,0.00000000;green,0.00000000;blue,0.00000000}, draw opacity = 1.0, rotate = 0.0},log basis y=10,ymajorgrids = true,ytick = {0.0031622776601683794,1.0,316.22776601683796,100000.0,3.162277660168379e7,1.0e10},yticklabels = {$10^{-2.5}$,$10^{0.0}$,$10^{2.5}$,$10^{5.0}$,$10^{7.5}$,$10^{10.0}$},ytick align = inside,yticklabel style = {font = {\fontsize{8 pt}{10.4 pt}\selectfont}, color = {rgb,1:red,0.00000000;green,0.00000000;blue,0.00000000}, draw opacity = 1.0, rotate = 0.0},y grid style = {color = {rgb,1:red,0.00000000;green,0.00000000;blue,0.00000000},
draw opacity = 0.1,
line width = 0.5,
solid},axis lines* = left,y axis line style = {color = {rgb,1:red,0.00000000;green,0.00000000;blue,0.00000000},
draw opacity = 1.0,
line width = 1,
solid},    xshift = 0.0mm,
    yshift = 0.0mm,
    axis background/.style={fill={rgb,1:red,1.00000000;green,1.00000000;blue,1.00000000}}
,legend style = {color = {rgb,1:red,0.00000000;green,0.00000000;blue,0.00000000},
draw opacity = 1.0,
line width = 1,
solid,fill = {rgb,1:red,1.00000000;green,1.00000000;blue,1.00000000},font = {\fontsize{8 pt}{10.4 pt}\selectfont}},colorbar style={title=}, ymin = {3.8018939632056124e-5}, width = \figurewidth]\addplot+ [color = {rgb,1:red,0.00000000;green,0.60560316;blue,0.97868012},
draw opacity = 1.0,
line width = 2,
solid,mark = *,
mark size = 2.0,
mark options = {
    color = {rgb,1:red,0.00000000;green,0.00000000;blue,0.00000000}, draw opacity = 1.0,
    fill = {rgb,1:red,0.00000000;green,0.60560316;blue,0.97868012}, fill opacity = 1.0,
    line width = 1,
    rotate = 0,
    solid
},forget plot]coordinates {
(10.0, 130.21121919632898)
(20.0, 49.3867665374658)
(30.0, 2.9775535829476745)
(40.0, 116.16275846012589)
(50.0, 2.49873801222045e8)
(60.0, 15.404719082991571)
(70.0, 7.142606734848604)
(80.0, 3.6696066455128067)
(90.0, 2.039062965809409)
(100.0, 1.2053401106286032)
};
\addplot+ [color = {rgb,1:red,0.88887350;green,0.43564919;blue,0.27812294},
draw opacity = 1.0,
line width = 2,
solid,mark = square*,
mark size = 2.0,
mark options = {
    color = {rgb,1:red,0.00000000;green,0.00000000;blue,0.00000000}, draw opacity = 1.0,
    fill = {rgb,1:red,0.88887350;green,0.43564919;blue,0.27812294}, fill opacity = 1.0,
    line width = 1,
    rotate = 0,
    solid
},forget plot]coordinates {
(10.0, 31.925562312991424)
(20.0, 15.576218219945867)
(30.0, 0.9247256366099946)
(40.0, 0.12433358184297717)
(50.0, 16.674261977347705)
(60.0, 6.708820344564353)
(70.0, 3.106855342892232)
(80.0, 1.5947547449980395)
(90.0, 0.8855342803950625)
(100.0, 0.5231748873209665)
};
\addplot+ [color = {rgb,1:red,0.24222430;green,0.64327509;blue,0.30444865},
draw opacity = 1.0,
line width = 2,
solid,mark = diamond*,
mark size = 2.0,
mark options = {
    color = {rgb,1:red,0.00000000;green,0.00000000;blue,0.00000000}, draw opacity = 1.0,
    fill = {rgb,1:red,0.24222430;green,0.64327509;blue,0.30444865}, fill opacity = 1.0,
    line width = 1,
    rotate = 0,
    solid
},forget plot]coordinates {
(10.0, 2.9417023632189734)
(20.0, 2.380391179415537)
(30.0, 0.1383545357589604)
(40.0, 0.018433794291070227)
(50.0, 0.0038635591995634317)
(60.0, 1.937529613820747)
(70.0, 0.8957421011521329)
(80.0, 0.45922106991531075)
(90.0, 0.25475950034399425)
(100.0, 0.15040332362729955)
};
\end{axis}

\end{tikzpicture}
\end{tiny}
}
\vspace{0.5cm}
\subfigure[$\protect \kappa_{M,N}^\epsilon$, equispaced points]{
\begin{tiny}
\begin{tikzpicture}[]
\begin{axis}[height = \figureheight, ylabel = {}, xmin = {7.300000000000001}, xmax = {102.7}, ymax = {2.4547089156850338e10}, ymode = {log}, xlabel = {N}, unbounded coords=jump,scaled x ticks = false,xlabel style = {font = {\fontsize{8 pt}{10.4 pt}\selectfont}, color = {rgb,1:red,0.00000000;green,0.00000000;blue,0.00000000}, draw opacity = 1.0, rotate = 0.0},xmajorgrids = true,xtick = {20.0,40.0,60.0,80.0,100.0},xticklabels = {$20$,$40$,$60$,$80$,$100$},xtick align = inside,xticklabel style = {font = {\fontsize{8 pt}{10.4 pt}\selectfont}, color = {rgb,1:red,0.00000000;green,0.00000000;blue,0.00000000}, draw opacity = 1.0, rotate = 0.0},x grid style = {color = {rgb,1:red,0.00000000;green,0.00000000;blue,0.00000000},
draw opacity = 0.1,
line width = 0.5,
solid},axis lines* = left,x axis line style = {color = {rgb,1:red,0.00000000;green,0.00000000;blue,0.00000000},
draw opacity = 1.0,
line width = 1,
solid},scaled y ticks = false,ylabel style = {font = {\fontsize{8 pt}{10.4 pt}\selectfont}, color = {rgb,1:red,0.00000000;green,0.00000000;blue,0.00000000}, draw opacity = 1.0, rotate = 0.0},log basis y=10,ymajorgrids = true,ytick = {0.01,1.0,100.0,10000.0,1.0e6,1.0e8,1.0e10},yticklabels = {$10^{-2}$,$10^{0}$,$10^{2}$,$10^{4}$,$10^{6}$,$10^{8}$,$10^{10}$},ytick align = inside,yticklabel style = {font = {\fontsize{8 pt}{10.4 pt}\selectfont}, color = {rgb,1:red,0.00000000;green,0.00000000;blue,0.00000000}, draw opacity = 1.0, rotate = 0.0},y grid style = {color = {rgb,1:red,0.00000000;green,0.00000000;blue,0.00000000},
draw opacity = 0.1,
line width = 0.5,
solid},axis lines* = left,y axis line style = {color = {rgb,1:red,0.00000000;green,0.00000000;blue,0.00000000},
draw opacity = 1.0,
line width = 1,
solid},    xshift = 0.0mm,
    yshift = 0.0mm,
    axis background/.style={fill={rgb,1:red,1.00000000;green,1.00000000;blue,1.00000000}}
,legend style = {color = {rgb,1:red,0.00000000;green,0.00000000;blue,0.00000000},
draw opacity = 1.0,
line width = 1,
solid,fill = {rgb,1:red,1.00000000;green,1.00000000;blue,1.00000000},font = {\fontsize{8 pt}{10.4 pt}\selectfont}},colorbar style={title=}, ymin = {0.0004073802778041126}, width = \figurewidth]\addplot+ [color = {rgb,1:red,0.00000000;green,0.60560316;blue,0.97868012},
draw opacity = 1.0,
line width = 2,
solid,mark = *,
mark size = 2.0,
mark options = {
    color = {rgb,1:red,0.00000000;green,0.00000000;blue,0.00000000}, draw opacity = 1.0,
    fill = {rgb,1:red,0.00000000;green,0.60560316;blue,0.97868012}, fill opacity = 1.0,
    line width = 1,
    rotate = 0,
    solid
},forget plot]coordinates {
(10.0, 25301.388140473904)
(20.0, 144711.52775834224)
(30.0, 140875.88427095397)
(40.0, 6.054375457912393e6)
(50.0, 887801.0060024441)
(60.0, 2.8850517812146306e7)
(70.0, 979654.578837018)
(80.0, 2.7668472420613185e7)
(90.0, 2.646524645860379e6)
(100.0, 6.621693055745844e7)
};
\addplot+ [color = {rgb,1:red,0.88887350;green,0.43564919;blue,0.27812294},
draw opacity = 1.0,
line width = 2,
solid,mark = square*,
mark size = 2.0,
mark options = {
    color = {rgb,1:red,0.00000000;green,0.00000000;blue,0.00000000}, draw opacity = 1.0,
    fill = {rgb,1:red,0.88887350;green,0.43564919;blue,0.27812294}, fill opacity = 1.0,
    line width = 1,
    rotate = 0,
    solid
},forget plot]coordinates {
(10.0, 2522.7607419822393)
(20.0, 8564.972662044158)
(30.0, 216610.6632972105)
(40.0, 74251.92967697686)
(50.0, 1.1085557796565895e6)
(60.0, 70498.95231341937)
(70.0, 844287.4153118087)
(80.0, 1.0303039734722009e7)
(90.0, 7.514584026357119e7)
(100.0, 1.9483412450143942e6)
};
\addplot+ [color = {rgb,1:red,0.24222430;green,0.64327509;blue,0.30444865},
draw opacity = 1.0,
line width = 2,
solid,mark = diamond*,
mark size = 2.0,
mark options = {
    color = {rgb,1:red,0.00000000;green,0.00000000;blue,0.00000000}, draw opacity = 1.0,
    fill = {rgb,1:red,0.24222430;green,0.64327509;blue,0.30444865}, fill opacity = 1.0,
    line width = 1,
    rotate = 0,
    solid
},forget plot]coordinates {
(10.0, 516.5477951304134)
(20.0, 1313.6740022574006)
(30.0, 860.4253282331048)
(40.0, 4552.716793716999)
(50.0, 22158.846737107644)
(60.0, 2572.6317620593404)
(70.0, 9619.583717033913)
(80.0, 35697.090637494584)
(90.0, 131703.0482281782)
(100.0, 483705.4198083835)
};
\end{axis}

\end{tikzpicture}
\end{tiny}
}
\subfigure[$\protect \lambda_{M,N}^\epsilon$, equispaced points]{
\begin{tiny}
\begin{tikzpicture}[]
\begin{axis}[height = \figureheight, ylabel = {}, xmin = {7.300000000000001}, xmax = {102.7}, ymax = {2.4547089156850338e10}, ymode = {log}, xlabel = {N}, unbounded coords=jump,scaled x ticks = false,xlabel style = {font = {\fontsize{8 pt}{10.4 pt}\selectfont}, color = {rgb,1:red,0.00000000;green,0.00000000;blue,0.00000000}, draw opacity = 1.0, rotate = 0.0},xmajorgrids = true,xtick = {20.0,40.0,60.0,80.0,100.0},xticklabels = {$20$,$40$,$60$,$80$,$100$},xtick align = inside,xticklabel style = {font = {\fontsize{8 pt}{10.4 pt}\selectfont}, color = {rgb,1:red,0.00000000;green,0.00000000;blue,0.00000000}, draw opacity = 1.0, rotate = 0.0},x grid style = {color = {rgb,1:red,0.00000000;green,0.00000000;blue,0.00000000},
draw opacity = 0.1,
line width = 0.5,
solid},axis lines* = left,x axis line style = {color = {rgb,1:red,0.00000000;green,0.00000000;blue,0.00000000},
draw opacity = 1.0,
line width = 1,
solid},scaled y ticks = false,ylabel style = {font = {\fontsize{8 pt}{10.4 pt}\selectfont}, color = {rgb,1:red,0.00000000;green,0.00000000;blue,0.00000000}, draw opacity = 1.0, rotate = 0.0},log basis y=10,ymajorgrids = true,ytick = {0.01,1.0,100.0,10000.0,1.0e6,1.0e8,1.0e10},yticklabels = {$10^{-2}$,$10^{0}$,$10^{2}$,$10^{4}$,$10^{6}$,$10^{8}$,$10^{10}$},ytick align = inside,yticklabel style = {font = {\fontsize{8 pt}{10.4 pt}\selectfont}, color = {rgb,1:red,0.00000000;green,0.00000000;blue,0.00000000}, draw opacity = 1.0, rotate = 0.0},y grid style = {color = {rgb,1:red,0.00000000;green,0.00000000;blue,0.00000000},
draw opacity = 0.1,
line width = 0.5,
solid},axis lines* = left,y axis line style = {color = {rgb,1:red,0.00000000;green,0.00000000;blue,0.00000000},
draw opacity = 1.0,
line width = 1,
solid},    xshift = 0.0mm,
    yshift = 0.0mm,
    axis background/.style={fill={rgb,1:red,1.00000000;green,1.00000000;blue,1.00000000}}
,legend style = {color = {rgb,1:red,0.00000000;green,0.00000000;blue,0.00000000},
draw opacity = 1.0,
line width = 1,
solid,fill = {rgb,1:red,1.00000000;green,1.00000000;blue,1.00000000},font = {\fontsize{8 pt}{10.4 pt}\selectfont}},colorbar style={title=}, ymin = {0.0004073802778041126}, width = \figurewidth]\addplot+ [color = {rgb,1:red,0.00000000;green,0.60560316;blue,0.97868012},
draw opacity = 1.0,
line width = 2,
solid,mark = *,
mark size = 2.0,
mark options = {
    color = {rgb,1:red,0.00000000;green,0.00000000;blue,0.00000000}, draw opacity = 1.0,
    fill = {rgb,1:red,0.00000000;green,0.60560316;blue,0.97868012}, fill opacity = 1.0,
    line width = 1,
    rotate = 0,
    solid
},forget plot]coordinates {
(10.0, 37218.21265649642)
(20.0, 56427.49348287582)
(30.0, 595027.1375537893)
(40.0, 343414.67444908136)
(50.0, 8.06463217332578e6)
(60.0, 7.172151670958242e6)
(70.0, 1.1138584374467996e8)
(80.0, 1.1018925195034869e8)
(90.0, 1.0929630919907996e8)
(100.0, 1.0856911380199239e8)
};
\addplot+ [color = {rgb,1:red,0.88887350;green,0.43564919;blue,0.27812294},
draw opacity = 1.0,
line width = 2,
solid,mark = square*,
mark size = 2.0,
mark options = {
    color = {rgb,1:red,0.00000000;green,0.00000000;blue,0.00000000}, draw opacity = 1.0,
    fill = {rgb,1:red,0.88887350;green,0.43564919;blue,0.27812294}, fill opacity = 1.0,
    line width = 1,
    rotate = 0,
    solid
},forget plot]coordinates {
(10.0, 22338.565927367305)
(20.0, 36029.78548479154)
(30.0, 9078.679065017623)
(40.0, 251138.04850234234)
(50.0, 163834.03689082246)
(60.0, 6.344726477179758e6)
(70.0, 5.75700075879243e6)
(80.0, 5.297478750313365e6)
(90.0, 1.032613706586497e8)
(100.0, 1.0670804206220502e8)
};
\addplot+ [color = {rgb,1:red,0.24222430;green,0.64327509;blue,0.30444865},
draw opacity = 1.0,
line width = 2,
solid,mark = diamond*,
mark size = 2.0,
mark options = {
    color = {rgb,1:red,0.00000000;green,0.00000000;blue,0.00000000}, draw opacity = 1.0,
    fill = {rgb,1:red,0.24222430;green,0.64327509;blue,0.30444865}, fill opacity = 1.0,
    line width = 1,
    rotate = 0,
    solid
},forget plot]coordinates {
(10.0, 120.98621353130477)
(20.0, 190.25626192487925)
(30.0, 3084.5186452922767)
(40.0, 1133.49313616904)
(50.0, 528.2198167964717)
(60.0, 54825.662254486684)
(70.0, 40605.18358574318)
(80.0, 31373.75643239411)
(90.0, 25029.756454565424)
(100.0, 20474.971946008605)
};
\end{axis}

\end{tikzpicture}
\end{tiny}
}
\caption{The same as Fig.\ \ref{fig:samplingrate} but with threshold $\epsilon = 1e^{-8}$. For Legendre nodes (top row), the intermediate peaks are higher, corresponding to the larger upper bounds in \eqref{kappa_lambda_global_indirect}. However, the values settle down for increasing $N$, with a clear benefit for larger oversampling factors. This is in agreement with the limit \eqref{kappa_lambda_limit_indirect}. For equispaced samples, linear oversampling is not sufficient and, as in the previous figure, the values seem to approach $1/\epsilon$. This again corresponds to \eqref{kappa_lambda_global_indirect}.}\label{fig:samplingrate2}
\end{center}
\end{figure}
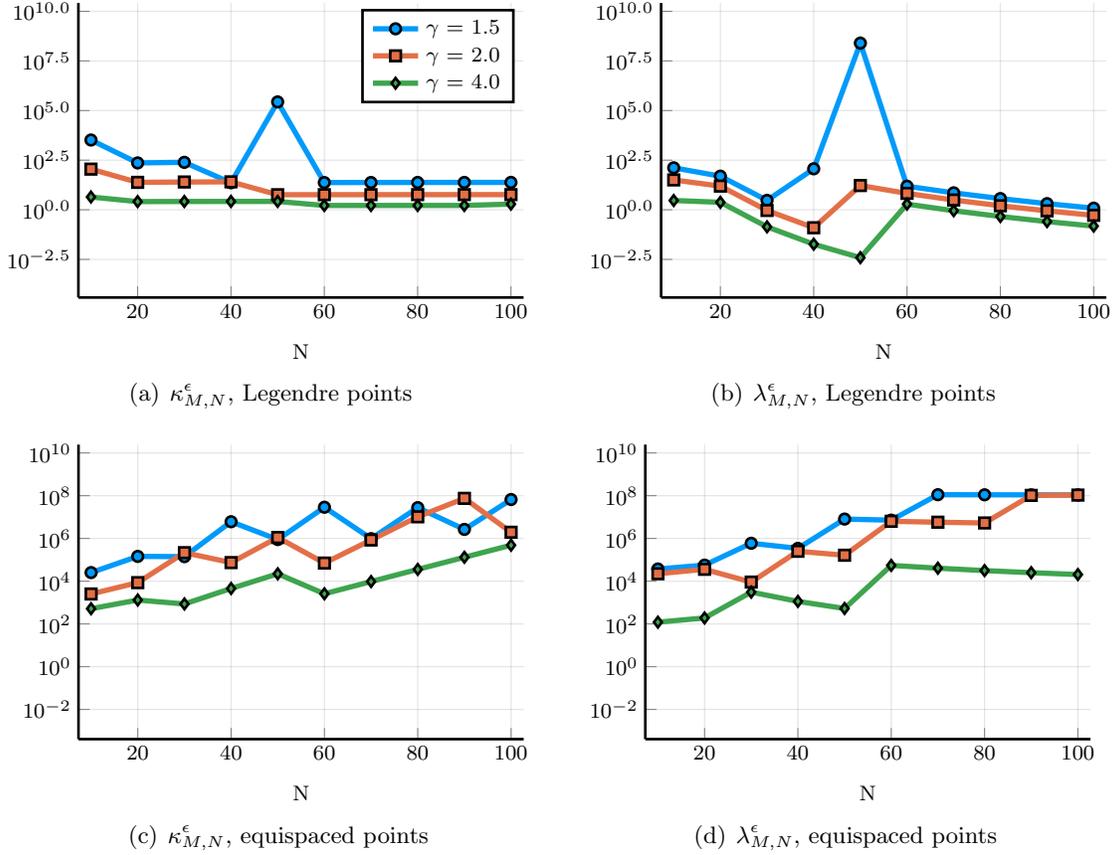

The results in Fig.\ \ref{fig:samplingrate} correspond to the threshold $\epsilon = 1e^{-5}$. For comparison, the experiment is repeated in Fig.\ \ref{fig:samplingrate2} for the smaller threshold $\epsilon = 1e^{-8}$. The latter figure illustrates the larger upper bound of \eqref{kappa_lambda_global_indirect}, on the order of $1/\epsilon$, in the pre-asymptotic regime. Still, for the case of Legendre nodes, linear oversampling is sufficient to reach the $\epsilon$-independent (and small) limit \eqref{kappa_lambda_limit_indirect}.

The constants were computed following the approach described in Remark \ref{r:SSRcompute}. We have run this experiment in higher precision arithmetic, in order to exclude the possibility of inaccuracies in their computation. An exponentially converging composite hp-graded quadrature rule was used to approximate the singular integrals that arise in the elements of the Gram matrix. Furthermore, in this case we also weighted the discrete samples in Legendre points by the square roots of the corresponding Gauss--Legendre quadrature weights: this discrete normalization ensures that $A' = B' = 1$ in \eqref{data_rich} and leads to slightly smaller values of the constants (and, correspondingly, smaller error in the approximation).

\section*{Acknowledgements}

A preliminary version of this work was presented during the Research Cluster on ``Computational Challenges in Sparse and Redundant Representations'' at ICERM in November 2014.  The authors would like to thank all the participants for the useful discussions and feedback received during the programme.  The first author would also like to thank Juan M.\ Cardenas and Sebastian Moraga.  The first author is supported by NSERC grant 611675, as well as an Alfred P.\ Sloan Research Fellowship  The second author is supported by FWO-Flanders projects G.0641.11, G.A004.14 and by KU Leuven project C14/15/055.

\bibliographystyle{abbrv}
\small
\bibliography{FramesStabilityRefs}

\begin{thebibliography}{10}

\bibitem{BAACHShannon}
B.~Adcock and A.~C. Hansen.
\newblock A generalized sampling theorem for stable reconstructions in
  arbitrary bases.
\newblock {\em J. Fourier Anal. Appl.}, 18(4):685--716, 2012.

\bibitem{BAACHAccRecov}
B.~Adcock and A.~C. Hansen.
\newblock Stable reconstructions in {H}ilbert spaces and the resolution of the
  {G}ibbs phenomenon.
\newblock {\em Appl. Comput. Harmon. Anal.}, 32(3):357--388, 2012.

\bibitem{BAACHOptimality}
B.~Adcock, A.~C. Hansen, and C.~Poon.
\newblock Beyond consistent reconstructions: optimality and sharp bounds for
  generalized sampling, and application to the uniform resampling problem.
\newblock {\em SIAM J. Math. Anal.}, 45(5):3114--3131, 2013.

\bibitem{FramesPart1}
B.~Adcock and D.~Huybrechs.
\newblock Frames and numerical approximation.
\newblock {\em SIAM Rev.}, 61(3):443--473, 2019.

\bibitem{FEStability}
B.~Adcock, D.~Huybrechs, and J.~Mart{\'\i}n-Vaquero.
\newblock On the numerical stability of {F}ourier extensions.
\newblock {\em Found. Comput. Math.}, 14(4):635--687, 2014.

\bibitem{boffi2015immersed}
D.~Boffi, N.~Cavallini, and L.~Gastaldi.
\newblock The finite element immersed boundary method with distributed
  {L}agrange multiplier.
\newblock {\em SIAM J. Numer. Anal.}, 53(6):2584--2604, 2015.

\bibitem{christensen2003introduction}
O.~Christensen.
\newblock {\em An Introduction to Frames and {R}iesz Bases}.
\newblock Applied and Numerical Harmonic Analysis. Birkh{\"a}user, 2nd edition,
  2016.

\bibitem{coppe2020splines}
V.~Copp{\'e} and D.~Huybrechs.
\newblock Efficient function approximation on general bounded domains using
  splines on a cartesian grid.
\newblock Technical Report arXiv:1911.07894, KU Leuven, 2020.

\bibitem{coppe2020az}
V.~Copp{\'e}, D.~Huybrechs, R.~Matthysen, and M.~Webb.
\newblock The {AZ} algorithm for least squares problems with a known incomplete
  generalized inverse.
\newblock {\em SIAM J. Mat. Anal. Appl.}, 2020.
\newblock To appear.

\bibitem{fix1973singular}
G.~J. Fix, S.~Gulati, and G.~I. Wakoff.
\newblock On the use of singular functions with finite elements approximations.
\newblock {\em J. Comput. Phys.}, pages 209--228, 1973.

\bibitem{golub1996matrix}
G.~H. Golub and C.~F. van Loan.
\newblock {\em Matrix computations}.
\newblock Johns Hopkins University Press, Baltimore, 3rd edition, 1996.

\bibitem{GrochenigMZineq}
K.~Gr\"{o}chenig.
\newblock Sampling, {M}arcinkiewicz--{Z}ygmund inequalities, approximation, and
  quadrature rules.
\newblock {\em arXiv:1909.07752}, 2019.

\bibitem{hansen2011}
A.~C. Hansen.
\newblock On the solvability complexity index, the n-pseudospectrum and
  approximations of spectra of operators.
\newblock {\em J. Amer. Math. Soc.}, 24(1):81--124, 2011.

\bibitem{Lindner2008}
E.~Heinemeyer, M.~Lindner, and R.~Potthast.
\newblock Convergence and numerics of a multisection method for scattering by
  three-dimensional rough surfaces.
\newblock {\em SIAM J. Numer. Anal.}, 46(4):1780--1798, 2008.

\bibitem{huybrechs2010fourier}
D.~Huybrechs.
\newblock On the {F}ourier extension of non-periodic functions.
\newblock {\em SIAM J. Numer. Anal.}, 47(6):4326--4355, 2010.

\bibitem{huybrechs2019wbm}
D.~Huybrechs and A.-E. Olteanu.
\newblock An oversampled collocation approach of the {W}ave {B}ased {M}ethod
  for {H}elmholtz problems.
\newblock {\em Wave Motion}, 87:92--105, 2019.

\bibitem{kasolis2015fictitious}
F.~Kasolis, E.~Wadbro, and M.~Berggren.
\newblock Analysis of fictitious domain approximations of hard scatterers.
\newblock {\em SIAM J. Numer. Anal.}, 2015(5):2347--2362, 2015.

\bibitem{lindner2006}
M.~Lindner.
\newblock {\em Infinite Matrices and their Finite Sections}.
\newblock Frontiers in Mathematics. Birkh\"auser Verlag, Basel, 2006.

\bibitem{LubinskyMarcinkiewicz}
D.~S. Lubinsky.
\newblock {\em Marcinkiewicz-Zygmund Inequalities: Methods and Results}, pages
  213--240.
\newblock Springer Netherlands, 1998.

\bibitem{LyonFast}
M.~Lyon.
\newblock A fast algorithm for {F}ourier continuation.
\newblock {\em {SIAM} J. Sci. Comput.}, 33(6):3241--3260, 2012.

\bibitem{matthysen2015fastfe}
R.~Matthysen and D.~Huybrechs.
\newblock Fast algorithms for the computation of {F}ourier extensions of
  arbitrary length.
\newblock {\em SIAM J. Sci. Comput.}, 38(2):A899--A922, 2016.

\bibitem{matthysen2017fastfe2d}
R.~Matthysen and D.~Huybrechs.
\newblock Function approximation on arbitrary domains using fourier frames.
\newblock {\em SIAM J. Numer. Anal.}, 56:1360--1385, 2018.

\bibitem{roache1978fourier}
P.~J. Roache.
\newblock A pseudo-spectral {FFT} technique for non-periodic problems.
\newblock {\em J. Comput. Phys.}, 27:204--220, 1978.

\bibitem{shirokoff2015volumepenalty}
D.~Shirokoff and J.-C. Nave.
\newblock A sharp-interface active penalty method for the incompressible
  {N}avier--{S}tokes equations.
\newblock {\em J. Sci. Comput.}, 62(1):53--77, 2015.

\end{thebibliography}

\end{document}